\documentclass[11pt]{article}
\usepackage{amssymb,amsmath,amsthm,mathtools}
\usepackage{graphicx}
\usepackage{booktabs,array,float}
\usepackage{enumitem}
\usepackage{xcolor}
\usepackage[hidelinks]{hyperref}
\allowdisplaybreaks

\usepackage{geometry}
\numberwithin{equation}{section}
\newtheorem{theorem}{Theorem}[section]
\newtheorem{lemma}[theorem]{Lemma}
\newtheorem{proposition}[theorem]{Proposition}
\newtheorem{corollary}[theorem]{Corollary}

\theoremstyle{remark}
\newtheorem{remark}[theorem]{Remark}
\theoremstyle{plain}

\def\Re{{\mathop{\rm Re}\,}}
\def\Im{{\mathop{\rm Im}\,}}
\newcommand{\R}{\mathbb R}
\newcommand{\C}{\mathbb C}
\newcommand{\N}{\mathcal N}
\newcommand{\E}{\mathcal E}
\newcommand{\Y}{\mathcal Y}
\newcommand{\X}{\mathcal X}
\newcommand{\Z}{\mathcal Z}
\newcommand{\LL}{\mathcal L}
\newcommand{\QQ}{\mathcal Q}
\newcommand{\Cd}{\mathfrak C}
\newcommand{\CP}{\mathrm{CP}}
\newcommand{\norm}[1]{\left\lVert #1\right\rVert}
\newcommand{\ip}[2]{\left(#1,#2\right)_\rho}
\newcommand{\dx}{\,dx}

\begin{document}
\title{Ordering error estimates under mesh-uniform regularity and optimal
two-mesh design for the KdV equation with boundary
dissipation\thanks{This work is supported by the National Natural Science
Foundation of China (grant No. 12171431).\medskip}}
\author{Kangsheng Liu\thanks{Academy for Advanced Interdisciplinary Studies,
Northeast Normal University, Changchun 130024, Jilin, China; School of
Mathematical Sciences, Zhejiang University, Hangzhou 310058, China
(email: ksliu@zju.edu.cn).}}
\date{}
\maketitle

\begin{center}
\begin{minipage}{6in}
{\bf Abstract.}
We study Galerkin error estimates and optimal two-mesh design for
the Korteweg--de Vries (KdV) equation on a bounded interval with boundary
dissipation. The ordering error is introduced to measure how much the
nonlinearity changes the discretization error. A weighted cubic Hermite
Galerkin scheme has mesh-uniform regularity: the discrete semigroups of
the linearized equation inherit Crandall--Pazy regularity and Kato
smoothing, with bounds uniform in the mesh size $h$. On a small ball of
the energy space, the discrete solution maps are also analytic in the
initial datum, with bounds uniform in $h$. On this ball the ordering error
is at most of order $h^{1/6}$ times the squared energy norm of the datum,
for all $t\ge0$. For $u_0$ on this ball and in the domain $D(A)$ of the linearized generator, the bound is of order
$h^2\norm{u_0}_{H^3}^{3/2}\norm{u_0}_{L^2}^{1/2}$.
\par\smallskip
A two-mesh computation adds a nonlinear correction on a coarse
mesh to a linear simulation on a fine mesh. Its error is the linear error
on the fine mesh plus the ordering error on the coarse mesh. Based on
these error estimates, we formulate a model of the optimal two-mesh
design: minimize the cost subject to an accuracy constraint given by the
error bounds. For data in a small ball of $L^2$, the problem is convex in logarithmic
mesh sizes. When the relative error bound for the coarsest allowed mesh exceeds the prescribed tolerance $\eta$,
the optimal pair and its KKT multipliers are unique. For smooth data
the problem is convex in the mesh sizes, and the optimal mesh ratio is
independent of the tolerance until the coarse mesh reaches its upper limit.
In both cases we characterize the small-amplitude limit of the relative
cost when the relative tolerance is proportional to the amplitude.
\vskip 3mm

{\bf Keywords.} Korteweg--de Vries equation, Hermite Galerkin method,
ordering error, mesh-uniform regularity, Crandall--Pazy semigroup,
two-mesh method, optimal mesh design
\vskip 4mm

{\bf 2020 Mathematics Subject Classification: } Primary 65M60, 65M15;
Secondary 35Q53, 47D06, 65M50, 65Y20, 90C46
\vskip 3mm
\end{minipage}
\end{center}

\section{Introduction}\label{sec:intro}

A small nonlinear wave is often simulated by a linear discrete model.
Two approximations are then made: the equation is linearized, and it is
discretized. They can be made in either order. Does the order matter?

We study this question for the Korteweg--de Vries (KdV)
equation~\cite{KdV} on an interval of fixed length $L>0$:
\begin{equation}\label{eq:pde}
 u_t+u_{xxx}+uu_x=0,\qquad
 u(0,t)=u(L,t)=u_x(L,t)=0,\qquad u(x,0)=u_0(x).
\end{equation}
We call $u_0$ the initial datum, or simply the datum. This paper makes
three contributions: a concept, a method of error estimation and an
application.

\paragraph{The ordering error.}
Our first contribution is a concept. Let $v$ solve the
linearized equation with the same datum, and let $u_h$ and $v_h$ be the
Galerkin approximations of $u$ and $v$. Linearizing first gives the path
$u\to v\to v_h$. Discretizing first gives the path $u\to u_h\to v_h$. Both
paths end at $v_h$, so the total error $u-v_h$ is the same. The order
decides how this error is divided between linearization and
discretization. Changing the order moves an amount from one share to the other,
\[
 \Cd_h=(u-v)-(u_h-v_h)=(u-u_h)-(v-v_h).
\]
We call $\Cd_h$ the \emph{ordering error}. It vanishes at $t=0$. Let
$\varepsilon$ be the size of the datum and $h^p$ a mesh factor. Separate
bounds for the four shares give only
$\Cd_h=O(\min\{\varepsilon^2,h^p\varepsilon\})$. Is the ordering error
of the product form $h^p\varepsilon^2$? Does this hold for all time?

\paragraph{Error estimates based on mesh-uniform regularity.}
Our second contribution is a method of error estimation. It rests on
\emph{mesh-uniform regularity}: regularity properties inherited by the
scheme with bounds uniform in the mesh size.
For a weighted Galerkin scheme, called the W-scheme below, the answer to
both questions is yes. Let $H=L^2((0,L),(L+x)\dx)$, and let
$\lambda=3\pi^2/(4L^3)$ be the energy decay rate described below. For
every $0<\nu<\lambda$ there are $r,C>0$, depending only on $L$, $\nu$ and
the mesh regularity, such that for real data
\begin{equation}\label{eq:headline}
 \sup_{t\ge0}e^{\nu t}\norm{\Cd_h(t)}_H\le Ch^{1/6}\norm{u_0}_H^2
 \qquad\text{if }\norm{u_0}_H\le r\text{ and }0<h\le\min\{1,L\}.
\end{equation}
The bound holds on a whole energy ball and for all time
(Theorem~\ref{thm:ordering}). No higher norm of the datum enters. The
full and linear errors cannot do this. The kernel of the initial
projection contains nonzero smooth functions. Hence $u-u_h$ and $v-v_h$
are not small at $t=0$ uniformly on a ball. Their difference is small
from $t=0$ on. For data in the domain $D(A)$ of the generator
$A=-\partial_x^3$ of the linearized equation, the linear and full errors
are of order $h^2$ on the whole half-line, and smallness is required only
in $H$. The ordering error then keeps a factor of the size of the datum:
it is at most $Ch^2\norm{Au_0}_H^{3/2}\norm{u_0}_H^{1/2}$
(Theorem~\ref{thm:DA}).

The starting point is the dissipation. We omit the transport term $u_x$,
which produces Rosier's critical lengths~\cite{Rosier}. With the
multiplier $L+x$, Rosier's multiplier method~\cite[Proposition~3.2]{Rosier}
shows that the linear flow on $H$ decays like $e^{-\lambda t}$ for every
$L$ and has Kato smoothing~\cite{Kato83} (Proposition~\ref{prop:rosier}).

Three properties enter, each uniformly in the mesh. Exponential stability
makes estimates valid for all time. Analytic dependence on the datum
separates the quadratic part of the solution from the rest. Smoothing in
time gives positive powers of $h$ for data of low regularity. The
delicate one is smoothing, because the KdV semigroup $T$ on a bounded
interval is not analytic. We use Crandall--Pazy regularity of index
$3/2$: $\norm{AT(t)}=O(t^{-3/2})$ as $t\to0$. A power bound with an
explicit index is what time-weighted error estimates require. The index
$3/2$ is attained: a wave packet $\chi e^{ikx}$ at time $t\sim k^{-2}$
gives $\norm{AT(t)}\ge ct^{-3/2}$. It is inherited, uniformly in the
mesh, by the cubic Hermite Galerkin scheme in the weighted space $H$
(Theorem~\ref{thm:cp}).

The energy identity gives a resolvent estimate of order $|\tau|^{-2/3}$,
uniform in the mesh. The Hilbert direct sum of the discrete generators
turns it into mesh-uniform regularity. By Wakaiki's
theorem~\cite{WakaikiCP}, a family of contraction semigroups on Hilbert
spaces is uniformly of Crandall--Pazy class if and only if it satisfies a
uniform resolvent decay estimate (Proposition~\ref{prop:uniformCP}). This
is the regularity counterpart of the characterization of uniform
exponential stability by a uniform resolvent bound, due to Liu and
Zheng~\cite{LiuZheng94}, which can be read as the
Gearhart--Huang--Pr\"uss theorem~\cite{Gearhart,Huang,Pruss} for the
direct sum. For the ordering error the key is a cancellation. In its
mild form the linear error $v-v_h$, which is not small uniformly on a
ball, drops out. Mixed time norms of type $L^4$ and $L^{4/3}$ control the
remaining terms on a whole energy ball.

The standard Galerkin method in $L^2(0,L)$ on the same trial space has
high-frequency spurious eigenmodes. Mesh-uniform regularity fails for
it, and so do the corresponding error estimates. The theoretical results,
the code of their Lean~4 check and the files of the computations are in
[arXiv:2611.12585v2].

\paragraph{Optimal two-mesh design.}
Our third contribution is an application: an optimization model for the
two-mesh computation, whose constraint is given by the error estimates.
A nonlinear simulation costs more than a linear one, since every time
step requires an iteration. It is therefore natural to compute the
nonlinear correction on a coarser mesh. Let $h\le\hat h$ be the sizes of
a fine and a coarse mesh. The two-mesh sum
$\tilde u=v_h+(u_{\hat h}-v_{\hat h})$ adds a nonlinear correction
computed on the coarse mesh to the linear simulation on the fine mesh. Its
error is $u-\tilde u=(v-v_h)+\Cd_{\hat h}$. The fine mesh serves the
linear part. The coarse mesh serves the ordering error, which carries the
factor $\varepsilon^2$. The three parts of the algebraic sum $\tilde u$ do not depend on
each other and can be computed in parallel.

The error estimates yield the explicit nonlinear program \eqref{eq:design}
for minimizing serial cost at a prescribed accuracy. In logarithmic mesh
sizes, the cost remains strictly convex and the feasible set becomes convex.
When the relative error bound for the coarsest allowed mesh exceeds the prescribed tolerance $\eta$ (condition~\eqref{eq:hbar}), the
Karush--Kuhn--Tucker (KKT) conditions~\cite[Section~4.2]{BSS} characterize
the unique optimal mesh pair (Theorem~\ref{prop:design}).
Corollary~\ref{cor:design-structure} describes how the two optimal mesh
sizes are related. Corollary~\ref{cor:design-small} gives the small-amplitude
limit of the relative cost when the relative tolerance is $\eta=\rho\varepsilon$ with fixed $\rho>0$.

For data in $D(A)$ both $v-v_h$ and $\Cd_h$ are of order $h^2$, and the design
problem is convex. The optimal solution is explicit
(Theorem~\ref{cor:design-smooth}). The optimal mesh ratio does not depend
on the accuracy unless the coarse mesh reaches its upper limit. For the same scaling
$\eta=\rho\varepsilon$, the relative cost has the same limit as for data
in $H$ (Theorem~\ref{cor:design-smooth}(c)).

The design uses two bounds only: a bound for the linear error on the fine
mesh and a bound of product form for the ordering error on the coarse
mesh. Wherever a scheme has such bounds, the same design applies. The
computations in Section~\ref{sec:comp} show how it performs for the KdV
equation.

\paragraph{Related work.}
We are not aware of an earlier study of the ordering error, nor of
mesh-uniform analyticity of solution maps for discretizations of KdV.
Analytic dependence on the data is classical for
KdV~\cite{Bourgain,KPV,BSZ}. On a bounded interval, Bona, Sun and
Zhang~\cite{BSZ} proved it on finite time intervals. Here the Taylor
remainder is uniform for all time, and it is obtained by a Catalan
majorant for the Duhamel power series. In parabolic problems,
exponential decay and an integrable smoothing kernel allow a contraction
in weighted trajectory spaces~\cite{Henry}. For KdV the derivative loss
requires Kato smoothing at the energy level.

Control and stabilization of KdV on a bounded interval have been studied
intensively~\cite{Rosier,PMVZ,RZ,RZsurvey,Cerpa}. Chu, Coron, Shang and
Tang~\cite{CCST} proved that the semigroup generated by
$-\partial_x-\partial_x^3$, with the boundary conditions of
\eqref{eq:pde}, is of Gevrey class $s$ for every $s>3/2$. Resolvent decay
of order $|\tau|^{-\beta}$ on the imaginary axis gives every Gevrey class
of order greater than $1/\beta$~\cite[Chapter~5]{Taylor89}. It also makes
the semigroup differentiable~\cite{CrandallPazy69}, with
$\norm{AT(t)}=O(t^{-\gamma})$ for every $\gamma>1/\beta$ on Banach
spaces~\cite[Theorem~2]{EEE95}. Wakaiki~\cite{WakaikiCP}
proved the exact characterization on Hilbert spaces: resolvent decay of
order $|\tau|^{-\beta}$ is equivalent to $\norm{AT(t)}=O(t^{-1/\beta})$
as $t\to0$.

Whether discretizations inherit the decay of a damped system has been
studied extensively. For the wave equation, numerical waves of high
frequency may escape boundary observation and damping, and the decay is
then not uniform in the mesh size. Many remedies have been proposed; see
the survey~\cite{Zuazua05} and the references therein. For the wave
equation with boundary feedback, Liu and Guo~\cite{LiuGuo} proposed an
order-reduction finite difference scheme that preserves the exponential
decay uniformly in the mesh size. For the
critical generalized KdV equation with localized damping, Pazoto,
Sep\'ulveda and Vera Villagr\'an~\cite{PSV} proved exponential decay for
small data, uniformly in the mesh parameters. Parabolic semigroups are
analytic, and their discretizations have a natural uniform regularity. For
non-analytic semigroups we are not aware of mesh-uniform Crandall--Pazy
regularity results.

Galerkin, Petrov--Galerkin, discontinuous Galerkin and finite-difference
schemes for KdV were analyzed in~\cite{BDK,DK,Shen,BCKX,HKR,CLR}. Most of
these works concern the periodic problem or the real line. Their error
estimates hold on finite time intervals and require a higher norm of the
solution, or of the datum, than the norm of the error. Our error bounds
hold on the whole half-line and decay exponentially. In
Theorems~\ref{thm:linear} and~\ref{thm:ordering}, the datum is measured in
the norm of the error.

Two-grid methods~\cite{Xu94,Xu96,MarionXu,DWW} combine a nonlinear
computation on a coarse mesh with simpler computations on a fine mesh.
For evolution problems their error estimates hold on bounded time
intervals, or on $[t_0,\infty)$ with constants depending on $t_0$ and on
the datum~\cite{MarionXu,DWW}.

\paragraph{Organization.}
Section~\ref{sec:continuous} treats the continuous problem.
Section~\ref{sec:discrete} introduces the W-scheme and proves its
mesh-uniform properties. The error estimates are given in
Section~\ref{sec:samenorm}. Section~\ref{sec:design} states and
solves the optimal design problem. Section~\ref{sec:comp} reports the
computations. Further estimates, including estimates in stronger norms,
are in [arXiv:2610.39269v2]. Table~\ref{tab:notation} collects the main notation.
The letter $C$ denotes a positive constant that may change from line to
line. Unless stated otherwise, it depends only on $L$, on the
mesh-regularity constant and on the decay rates.

\begin{table}[H]
\centering\small
\caption{Notation.}\label{tab:notation}
\begin{tabular}{@{}l>{\raggedright\arraybackslash}p{0.56\textwidth}l@{}}
\toprule
Symbol & Meaning & Defined in\\
\midrule
$L$; $h$, $N$ & length of the interval $(0,L)$; mesh size, and number of elements of a uniform mesh, $h=L/N$ & Section~\ref{sec:scheme}\\
$H$; $\norm\cdot_2$, $|\cdot|_1$ & weighted space $L^2((0,L),(L+x)\,dx)$; norm of $L^2(0,L)$, $|f|_1=\norm{f'}_2$ & Section~\ref{sec:setting}\\
$A$, $B$, $T(t)$; $\lambda$ & generator $-\partial_x^3$, $B=-A$, semigroup of the linearized problem; energy decay rate & Section~\ref{sec:setting}\\
$V_\theta$ & fractional scale $D(B^{\theta/3})$ & \eqref{eq:scale}\\
$\E_\nu$; $\X_0$ & energy norm of a trajectory; energy trajectory space & \eqref{eq:Enorm}, \eqref{eq:X0}\\
$X_h$; $A_h$, $T_h(t)$, $\Pi_h$ & trial space; W-scheme: weighted Galerkin operator, its semigroup, $H$-orthogonal projection & Section~\ref{sec:scheme}\\
$u$, $v$; $u_h$, $v_h$ & nonlinear and linear solutions: continuous, W-scheme & Section~\ref{sec:orders}\\
$\Cd_h$; $L_h(t)$ & ordering error; $T_h(t)\Pi_h-T(t)$ & Sections~\ref{sec:intro}, \ref{sec:orders}\\
$D_h$ & discretization error of the quadratic Taylor coefficient & \eqref{eq:Dh}\\
$\hat h$; $\tilde u$ & coarse mesh size; two-mesh sum $v_h+(u_{\hat h}-v_{\hat h})$ & Section~\ref{sec:design}\\
$\Phi$, $\Psi$ & error terms of the linear part and of the correction & \eqref{eq:PhiPsi}, \eqref{eq:smoothPhiPsi}\\
$J_{\rm L}$, $J_{\rm N}$; $d_{\rm L}$, $d_{\rm N}$ & computing time of a linear and of a nonlinear run; their time per unknown and time step & Section~\ref{sec:design-cost}\\
$M$; $\gamma$ & bound $\norm{Au_0}_H\le M\norm{u_0}_H$ for data in $D(A)$; optimal mesh ratio for such data when $h^*<\hat h^*<\bar h$ & Section~\ref{sec:design-smooth}\\
$x_\wedge$ & $\min\{1,x\}$ for $x\ge0$, for instance $t_\wedge$ and $(t-s)_\wedge$ & throughout\\
$\langle\omega\rangle$ & $(1+\omega^2)^{1/2}$ & throughout\\
\bottomrule
\end{tabular}
\end{table}

\section{The continuous problem}\label{sec:continuous}
\setcounter{equation}{0}

In this section, we treat the continuous problem. We first introduce the
energy space, the adjoint and the fractional scale. We then state that
the solution map is analytic, and we prove Crandall--Pazy regularity of
exact index $3/2$. The proof of
Theorem~\ref{thm:ctcp}(i) uses the discrete resolvent bound, which is
proved independently in Section~\ref{sec:uniformcp}. The last subsection
proves the analyticity.

\subsection{Energy, adjoint and fractional spaces}\label{sec:setting}
Set $\rho(x)=L+x$ and $H=L^2((0,L),\rho\dx)$, with inner product
$\ip fg=\int_0^Lf\bar g\rho\dx$ and norm $\norm\cdot_H$. We write
$\norm\cdot_2$ for the unweighted $L^2$ norm and $|f|_1=\norm{f'}_2$. A
prime always denotes a spatial derivative. Define
\begin{equation}\label{eq:operator}
 A=-\partial_x^3,\qquad
 D(A)=\{f\in H^3(0,L): f(0)=f(L)=f'(L)=0\},
\end{equation}
and $F(f,g)=-\frac12(fg)'$, $\N(f)=F(f,f)$. Then \eqref{eq:pde} reads
$u_t=Au+\N(u)$. The bilinear map $F$ involves no complex conjugation. Complex
spaces are used in the analyticity arguments. The physical solutions are
real.

For $f\in H^2$ with $f(0)=f(L)=f'(L)=0$ and $\psi\in H_0^1(0,L)$ put
\begin{equation}\label{eq:form}
 b(f,\psi)=-\int_0^Lf''\,\overline{(\rho\psi)'}\dx .
\end{equation}
Integrating by parts once, we obtain $\ip{Af}\psi=-b(f,\psi)$ for $f\in D(A)$.
Put
\[
 a=\frac32,\qquad p=\frac{\sqrt{2L}\,L}{\pi},\qquad
 \lambda=\frac a{p^2}=\frac{3\pi^2}{4L^3}.
\]
Two further integrations by parts yield
\begin{equation}\label{eq:energy}
 \Re b(f,f)=a|f|_1^2+\frac L2|f'(0)|^2,\qquad
 \norm f_H\le p|f|_1 .
\end{equation}
The inequality combines $\rho\le2L$ with Poincar\'e's inequality.

The next proposition is a modification of~\cite[Propositions~3.1 and~3.2]{Rosier}.
There the space is $L^2(0,L)$, the equation contains the transport term,
and the multiplier estimates hold on finite time intervals.

\begin{proposition}[Contraction, decay and Kato smoothing]\label{prop:rosier}
The operator $A$ is invertible, $-A$ is maximal accretive, and $A$
generates a contraction semigroup $T(t)=e^{tA}$ on $H$. For $f\in H$ and
$t\ge0$,
\begin{equation}\label{eq:ctenergy}
 \norm{T(t)}_{H\to H}\le e^{-\lambda t},\qquad
 \norm{T(t)f}_H^2+2a\int_0^t|T(s)f|_1^2\,ds\le\norm f_H^2 .
\end{equation}
\end{proposition}
\begin{proof}
By \eqref{eq:energy}, $-\Re\ip{Af}f=\Re b(f,f)\ge0$ for $f\in D(A)$. For
$s\ge0$ the homogeneous third-order boundary value problem $(s-A)f=0$
therefore has only the zero solution. Hence its three-by-three boundary matrix
is invertible. Variation of constants then solves $(s-A)f=g$ for every
$g\in H$, with $\norm f_{H^3}\le C\norm g_H$. Hence $A$ is invertible,
$-A$ is maximal accretive, and $A$ generates a contraction semigroup.

Let $f\in D(A)$ and $u(t)=T(t)f$. By \eqref{eq:energy},
\[
 \frac d{dt}\norm{u(t)}_H^2=-2\Re b(u,u)=-2a|u|_1^2-L|u'(0)|^2 .
\]
This is the multiplier identity~\cite[(3.16)]{Rosier} with $q=L+x$,
without the transport term. By the Poincar\'e inequality in \eqref{eq:energy} and $\lambda=a/p^2$, we have $2a|u|_1^2\ge2\lambda\norm u_H^2$. Gronwall's inequality then yields the first estimate. Integrating the identity over $(0,t)$, we obtain the second. Both extend
to $f\in H$ by density.
\end{proof}
The second estimate is the Kato smoothing of the linear flow.
Compare~\cite[(3.18)]{Rosier}.

Three integrations by parts identify the $H$-adjoint of $A$:
\begin{equation}\label{eq:adjoint}
\begin{gathered}
 A^*g=\rho^{-1}(\rho g)''',\qquad
 D(A^*)=\{g\in H^3:g(0)=g(L)=g'(0)=0\},\\
 -\Re\ip{A^*g}g=a|g|_1^2+L|g'(L)|^2,\qquad g\in D(A^*).
\end{gathered}
\end{equation}

Define the fractional scale
\begin{equation}\label{eq:scale}
 B=-A,\qquad V_\theta=D(B^{\theta/3}),\qquad
 \norm f_{V_\theta}=\norm{B^{\theta/3}f}_H,\qquad 0\le\theta\le3 .
\end{equation}
In particular $V_0=H$, $V_1=H_0^1(0,L)$,
\[
 V_2=\{f\in H^2:f(0)=f(L)=f'(L)=0\},\qquad V_3=D(A),
\]
and at these integer indices the norms are equivalent to the Sobolev
norms. By the half-plane functional calculus for Hilbert-space contraction semigroups, we have $\norm{B^{is}}\le e^{\pi|s|/2}$ for real
$s$~\cite[Theorem~11.12]{Haase}. By interpolation of fractional powers and of boundary spaces~\cite{Seeley}, we obtain
\begin{equation}\label{eq:interp}
 [V_r,V_s]_\alpha=V_{(1-\alpha)r+\alpha s},\qquad
 0\le r,s\le3,\quad 0\le\alpha\le1 .
\end{equation}
Away from $\theta=1/2,3/2$, $V_\theta$ is the subspace of $H^\theta$
carrying exactly the boundary conditions of order $j<\theta-1/2$. At the
critical indices, $V_\theta$ is given by interpolation. For instance,
$V_{1/2}=H_{00}^{1/2}$. In all cases $V_\theta\hookrightarrow H^\theta$.
We shall often use the moment inequality
\begin{equation}\label{eq:moment}
 \norm f_{V_\theta}\le C\norm f_H^{1-\theta/s}\norm f_{V_s}^{\theta/s},
 \qquad 0\le\theta\le s\le3,\quad s>0 .
\end{equation}
For $1\le\theta\le3$ the nonlinearity satisfies
\begin{equation}\label{eq:product}
 \norm{F(f,g)}_{V_{\theta-1}}
 \le C_\theta\norm f_{V_\theta}\norm g_{V_\theta}.
\end{equation}
At integer indices this follows from the one-dimensional Sobolev product estimates and
the boundary conditions. For $\theta=2$, $(fg)'$ vanishes at both ends.
For $\theta=3$ we also have
$(fg)''(L)=f''(L)g(L)+2f'(L)g'(L)+f(L)g''(L)=0$. The estimates at the intermediate indices, including the critical ones, follow by bilinear complex interpolation.

\subsection{The analytic solution map}\label{sec:analytic}
Let $V$ be a Banach space and $\nu>0$. We define
\begin{gather*}
 \Y_\nu(V)=\Big\{w\in C([0,\infty);V):\
 \norm w_{\Y_\nu(V)}:=\sup_{t\ge0}e^{\nu t}\norm{w(t)}_V<\infty\Big\},\\
 C_\nu^0(V)=\Big\{w\in\Y_\nu(V):\
 \lim_{t\to\infty}e^{\nu t}\norm{w(t)}_V=0\Big\}.
\end{gather*}
The space $C_\nu^0(V)$ is a closed subspace of $\Y_\nu(V)$. For
$\sigma>0$ put
\begin{equation}\label{eq:Enorm}
 \E_\sigma(w)=\sup_{t\ge0}e^{\sigma t}\norm{w(t)}_H+
 \Big(\int_0^\infty e^{2\sigma t}|w(t)|_1^2\,dt\Big)^{1/2}.
\end{equation}
For fixed $\nu$, the energy trajectory space is
\begin{equation}\label{eq:X0}
 \X_0=C_\nu^0(H)\cap L^2(e^{2\nu t}dt;H_0^1),\qquad
 \norm w_{\X_0}=\E_\nu(w).
\end{equation}

\begin{theorem}[Analytic solution map]\label{thm:continuous}
Fix $0<\nu<\lambda$. There are $r>0$, $M>0$, both
depending only on $L$ and $\nu$, and continuous homogeneous polynomials
$P_n:H\to\X_0$ of degree $n$, with $P_1(f)=T(\cdot)f$, such
that
\begin{equation}\label{eq:Pn}
 \norm{P_n(f)}_{\X_0}\le M\Big(\frac{\norm f_H}{2r}\Big)^n,
 \qquad f\in H,\ n\ge1 .
\end{equation}
For complex $f\in H$ with $\norm f_H<r$ the series
\[
 \mathcal U(f)=\sum_{n\ge1}P_n(f)
\]
converges absolutely in the norm of $\X_0$, uniformly on this ball.
Its sum is the unique solution of the mild form of \eqref{eq:pde} in a
ball of $\X_0$ whose radius depends only on $L$ and $\nu$. It is
real for real $f$. For every integer $m\ge1$ there is $C_m$ such that
\begin{equation}\label{eq:taylor}
 \norm{\mathcal U(f)-\sum_{n=1}^mP_n(f)}_{\X_0}
 \le C_m\norm f_H^{m+1},\qquad \norm f_H<r .
\end{equation}
In particular, with $\varepsilon=\norm{u_0}_H<r$ and
$u=\mathcal U(u_0)$,
\begin{equation}\label{eq:continuous}
 \norm{u(t)}_H\le C\varepsilon e^{-\nu t},\qquad
 \norm{u(t)-T(t)u_0}_H\le C\varepsilon^2e^{-\nu t},\qquad
 t\ge0 .
\end{equation}
\end{theorem}

Relative to the norm of the datum, the linearization error is
$O(\varepsilon)$ for all time. The series also converges in the
weighted $L^2(H_0^1)$ norm. This form of the statement survives
discretization (Section~\ref{sec:uniformanalytic}).
Theorem~\ref{thm:continuous} is proved in Section~\ref{sec:analyticproof}.

\subsection{Crandall--Pazy regularity}\label{sec:cp}
Let $S$ be a $C_0$-semigroup with generator $G$ on a Hilbert space. It is
\emph{immediately differentiable} if $S(t)$ maps into $D(G)$ for every
$t>0$. For $0<\beta\le1$ it belongs to the \emph{Crandall--Pazy class}
$\CP_\beta$ if, in addition,
$\limsup_{t\to0+}t^{1/\beta}\norm{GS(t)}<\infty$. We then say that $S$
has \emph{Crandall--Pazy regularity} of index $1/\beta$.

\begin{theorem}[Wakaiki~{\cite[Theorem~2.3]{WakaikiCP}}]\label{thm:wakaiki}
Let $G$ generate a $C_0$-semigroup $S$ on a Hilbert space and let
$0<\beta\le1$. Then $S\in\CP_\beta$ if and only if there are
$\Omega,C>0$ such that $i\omega\in\rho(G)$ and
$\norm{R(i\omega,G)}\le C|\omega|^{-\beta}$ for all $|\omega|\ge\Omega$.
\end{theorem}
In~\cite{WakaikiCP} the generator is denoted by $-A$ and the semigroup by $e^{-tA}$. With $G=-A$ the
statement there is the one above. See also~\cite[Theorem~5.1]{Wakaiki} for
general resolvent growth functions.

\begin{theorem}[Crandall--Pazy regularity of exact index]\label{thm:ctcp}
The following hold with constants depending only on $L$.
\begin{enumerate}[label=\rm(\roman*)]
\item $i\R\subset\rho(A)$ and $\norm{R(i\tau,A)}\le C\langle\tau\rangle^{-2/3}$
 for $\tau\in\R$.
\item For $0\le\alpha\le1$ and $t>0$,
 $\norm{B^\alpha T(t)}_{H\to H}\le Ct_\wedge^{-3\alpha/2}e^{-\lambda t}$.
\item There are $c_1,t_0>0$ such that $\norm{AT(t)}_{H\to H}\ge c_1t^{-3/2}$
 for $0<t\le t_0$.
\end{enumerate}
In particular $T\in\CP_{2/3}$, and the index $3/2$ cannot be improved.
\end{theorem}
\begin{proof}
(i) By \eqref{eq:ctenergy}, every $z$ with $\Re z>-\lambda$ lies in
$\rho(A)$, and $\norm{R(i\tau,A)}\le\lambda^{-1}$. The decay bound in (i) follows from
the mesh-uniform bound \eqref{eq:resolvent} on uniform meshes
($\kappa=1$), since $R(i\tau,A_h)\Pi_hf\to R(i\tau,A)f$ for every
$f\in H$. The discrete resolvent bound \eqref{eq:resolvent} is proved in
Section~\ref{sec:uniformcp} without using the present theorem.

To prove the convergence, let $y=R(i\tau,A)f$. Then $y\in D(A)$ and
$\norm y_{H^3}\le C\norm{Ay}_H\le C(1+|\tau|)\norm f_H$. With
the residual $d_h$ of Lemma~\ref{lem:projection},
$(i\tau-A_h)\Pi_hy=\Pi_hf+d_h(y)$, so
$R(i\tau,A_h)\Pi_hf-\Pi_hy=-R(i\tau,A_h)d_h(y)$. For $r,\psi\in X_h$ put
$x=R(-i\tau,A_h^*)\psi$. Then $\ip{R(i\tau,A_h)r}\psi=\ip rx$. The real
part of $\ip{(-i\tau-A_h^*)x}x=\ip\psi x$ and \eqref{eq:energy} yield
$a|x|_1^2\le\norm\psi_H\norm x_H\le\lambda^{-1}\norm\psi_H^2$. Hence
$\norm{R(i\tau,A_h)r}_H\le(a\lambda)^{-1/2}\norm r_{-1,h}$. By
\eqref{eq:residual} with $s=3$,
\[
 \norm{R(i\tau,A_h)\Pi_hf-\Pi_hy}_H
 \le Ch\norm y_{H^3}\le Ch(1+|\tau|)\norm f_H\longrightarrow0 .
\]
Since $\Pi_hy\to y$, the convergence follows. Thus
$\norm{R(i\tau,A)f}_H=\lim_{h\to0}\norm{R(i\tau,A_h)\Pi_hf}_H
\le C\langle\tau\rangle^{-2/3}\norm f_H$, which is (i).

(ii) By (i) and Theorem~\ref{thm:wakaiki}, $\norm{BT(t)}\le Ct^{-3/2}$ on
some interval $(0,t_1]$. For $t_1\le t\le1$ we use
$BT(t)=T(t-t_1)BT(t_1)$. For $t\ge1$ we use $BT(t)=T(t-1)BT(1)$ and
\eqref{eq:ctenergy}. This proves (ii) for $\alpha=1$. The general case follows from the moment inequality $\norm{B^\alpha y}\le C\norm y^{1-\alpha}\norm{By}^\alpha$, applied to $y=T(t)x$.

(iii) Choose $\chi\in C_c^\infty(0,L)$ with $\norm\chi_H=1$ and set
$f_k=\chi e^{ikx}$, for real $k\ge1$. Then $f_k\in D(A)$ and
\[
 r_k:=(A-ik^3)f_k=e^{ikx}(-\chi'''-3ik\chi''+3k^2\chi'),\qquad
 \norm{r_k}_H\le C_0k^2 .
\]
By the contraction property and the variation of constants formula, we obtain
$\norm{T(t)f_k-e^{ik^3t}f_k}_H\le C_0tk^2$. Fix $b_0>0$ with
$C_0b_0\le1/4$ and take $t=b_0/k^2$. Since $A$ commutes with $T(t)$ on $D(A)$, we obtain
\[
 \norm{AT(t)f_k}_H\ge k^3\norm{T(t)f_k}_H-\norm{T(t)r_k}_H
 \ge\tfrac34k^3-C_0k^2\ge\tfrac12k^3
\]
for sufficiently large $k$. Since $k\ge1$ is arbitrary,
the choice $k=(b_0/t)^{1/2}$ gives the lower bound for every sufficiently small
$t$. This proves (iii).
\end{proof}

The proof of \eqref{eq:resolvent} in Section~\ref{sec:uniformcp} uses the
following construction. Let $x\in H_0^1$ with $\norm x_H=1$ and
$k=|x|_1$, so that $pk\ge1$. Fix $\gamma\in(0,L/(8p)]$ and put
$\ell=\gamma/k\le L/8$. We extend $x$ by zero, mollify at scale $\ell$
and cut off in end strips of that scale. This produces
$z\in C_c^\infty(0,L)$ with
\begin{equation}\label{eq:mollify}
 \norm{x-z}_2\le C\ell k,\qquad \norm{z'}_2\le Ck,\qquad
 \norm{z''}_2\le C\ell^{-1}k .
\end{equation}
The cutoff terms are controlled by
$\norm x_{L^2(0,c\ell)}\le c\ell\norm{x'}_{L^2(0,c\ell)}$ and its
analogue at the right end. Each derivative of the cutoff costs one factor
$\ell^{-1}$, which this strip estimate offsets.

Let $C_R$ be the
constant in Theorem~\ref{thm:ctcp}(i). A Neumann series around $i\omega$ shows that
\begin{equation}\label{eq:strip}
 z=-\sigma+i\omega\in\rho(A),\qquad
 \norm{R(z,A)}\le2C_R\langle\omega\rangle^{-2/3},\qquad
 0\le\sigma\le(2C_R)^{-1}\langle\omega\rangle^{2/3}.
\end{equation}

\subsection{Proof of Theorem~\ref{thm:continuous}}\label{sec:analyticproof}
The proof is divided into two steps. We first close the estimates at the
energy level in a trajectory space, and then sum the power series.

\emph{The energy level.} The pointwise map $F:L^2\times L^2\to H^{-1}$
is not bounded. For $f_n(x)=n^{1/2}\chi(n(x-x_*))$, with $\chi$ supported
inside the interval, $\norm{f_n}_2$ is bounded whereas
$\norm{(f_n^2)'}_{H^{-1}}\asymp n^{1/2}$. The estimates must therefore be
closed in a trajectory space. Let $V_1^*$ be the dual of $V_1$ with respect to the
weighted pairing. In all weak-source estimates, we equip $H_0^1(0,L)$ with
the equivalent norm $|\cdot|_1$ and use the corresponding dual norm,
denoted by $\norm\cdot_{V_1^*}$. For $y_t=Ay+g$ the weighted energy estimate yields
\begin{equation}\label{eq:forced}
 \E_\nu(y)\le C_{L,\nu}
 \big(\norm{y(0)}_H+\norm{e^{\nu t}g}_{L^2_t(V_1^*)}\big),
 \qquad 0<\nu<\lambda .
\end{equation}
To see this, we multiply the energy identity for $y$, which follows from \eqref{eq:energy}, by $e^{2\nu t}$ and absorb the
source into $(a-\nu p^2)|y|_1^2$, whose coefficient is positive. Smooth
approximation defines the forced solution for such weak sources. For
sources that are smooth and compactly supported in time, the solution is
classical; after the support of the source only the free flow remains,
which decays at the rate $\lambda>\nu$. Such sources are dense in the
weighted space $L^2_t(V_1^*)$. By \eqref{eq:forced}, the weighted
trajectories converge uniformly, and the limit keeps continuity and decay
at infinity.

For $w,\tilde w\in\X_0$, the weak product bound is
$\norm{F(w,\tilde w)}_{V_1^*}\le C\norm{w\tilde w}_2$. By Agmon's inequality, $\norm{w\tilde w}_2^2\le C\norm w_H|w|_1\norm{\tilde w}_H^2$. If
$M_w=\sup_te^{\nu t}\norm{w(t)}_H$, then
\begin{align*}
 \norm{e^{\nu t}w\tilde w}_{L^2_tL^2_x}^2
 &\le CM_wM_{\tilde w}^2\int_0^\infty e^{-\nu t}|w|_1\,dt
  \le C_{L,\nu}\E_\nu(w)^2\E_\nu(\tilde w)^2 ,
\end{align*}
by the Cauchy--Schwarz inequality with the factor $e^{-2\nu t}$. Thus the
maps
\begin{equation}\label{eq:lowbilinear}
 (\LL f)(t)=T(t)f,\qquad
 \QQ(w,\tilde w)(t)=\int_0^tT(t-s)F(w(s),\tilde w(s))\,ds
\end{equation}
are bounded as $\LL:H\to\X_0$ and $\QQ:\X_0\times\X_0\to\X_0$, over the
complex field as well as the real field. On finite time intervals $F(w,\tilde w)$ also belongs to $L^1_tH$. Hence the weak definition of $\QQ$ agrees with the
usual mild formula.

\emph{The power series.} Let $l$ and $q$ be the bounds for $\norm\LL$
and $\norm\QQ$ obtained above. They depend only on $L$ and $\nu$. The
recursion
\[
 P_1(f)=\LL f,\qquad P_n(f)=\sum_{j+k=n}\QQ(P_j(f),P_k(f))
\]
produces continuous homogeneous polynomials. If $C_n$ are the Catalan
numbers, then
\[
 \norm{P_n(f)}_{\X_0}\le C_{n-1}q^{n-1}l^n\norm f_H^n
 \le4^{n-1}q^{n-1}l^n\norm f_H^n .
\]
Put $r=(8ql)^{-1}$. Then $4^{n-1}q^{n-1}l^n=(4q)^{-1}(2r)^{-n}$, which is
\eqref{eq:Pn} with $M=(4q)^{-1}$. The series converges normally on the
ball $\norm f_H<r$. Rearranging the absolutely convergent
double series $\sum_{j,k}\QQ(P_j,P_k)$ shows that the sum solves
$u=\LL f+\QQ(u,u)$. Its norm is at most
$l\norm f/(1-4ql\norm f)\le2l\norm f\le(4q)^{-1}$. On the ball of radius
$(4q)^{-1}$ in $\X_0$, the map $w\mapsto\LL f+\QQ(w,w)$ is a
contraction. Indeed, for $w,\tilde w$ in that ball,
$\norm{\QQ(w,w)-\QQ(\tilde w,\tilde w)}\le q(\norm w+\norm{\tilde w})
\norm{w-\tilde w}\le\frac12\norm{w-\tilde w}$. Hence the sum is
the unique fixed point in that ball. Summing the geometric majorant over the degrees $n>m$, we obtain \eqref{eq:taylor}. By conjugation and uniqueness, the values are real for real data. No implicit function theorem is used. This completes the proof of Theorem~\ref{thm:continuous}.
\hfill$\square$

\section{The W-scheme and its mesh-uniform properties}\label{sec:discrete}
\setcounter{equation}{0}

In this section, we introduce the W-scheme and prove its
mesh-uniform properties. These are exponential stability, Kato smoothing,
Crandall--Pazy regularity and analyticity of the solution map.

Let a family of discrete semigroups $T_h$ on subspaces $X_h\subset H$,
with generators $A_h$, approximate a semigroup $T$. The family has
\emph{mesh-uniform regularity} if the semigroups $T_h$ inherit a
regularity property of $T$, with bounds independent of $h$. Two kinds
matter here: Crandall--Pazy regularity of index $3/2$,
$\sup_h\sup_{0<t\le1}t^{3/2}\norm{A_hT_h(t)}<\infty$, and Kato smoothing,
$\int_0^\infty|T_h(t)x|_1^2\,dt\le C\norm x_H^2$ for $x\in X_h$, with $C$
independent of $h$. The W-scheme has both
(Theorem~\ref{thm:cp} and \eqref{eq:disKato}). The error estimates of
Section~\ref{sec:samenorm} rest on them.

\subsection{The scheme and the projection}\label{sec:scheme}
Let $0=x_0<\cdots<x_N=L$ be a mesh with
$h=\max_j(x_j-x_{j-1})\le\min\{1,L\}$ and
$h/\min_j(x_j-x_{j-1})\le\kappa$. Let $X_h\subset V_2$ be the space of
$C^1$ piecewise cubic polynomials that satisfy the three boundary
conditions of \eqref{eq:pde}, and let $I_h$ be the cubic Hermite
interpolant~\cite[Chapter~2]{Ciarlet}. Let $\Pi_h$ be the $H$-orthogonal
projection onto $X_h$. The
weighted Galerkin operator $A_h$ on $X_h$ is defined by
\begin{equation}\label{eq:scheme}
 \ip{A_hw}\psi=-b(w,\psi),\qquad w,\psi\in X_h,
\end{equation}
with the form \eqref{eq:form}. Put $B_h=-A_h$, $T_h(t)=e^{tA_h}$ and
$F_h=\Pi_hF$. The W-scheme is
\begin{equation}\label{eq:Wscheme}
 u_{h,t}=A_hu_h+F_h(u_h,u_h),\qquad u_h(0)=\Pi_hu_0 .
\end{equation}
All results concern the semidiscretization in space. Time remains continuous. We do
not change the initialization and do not add artificial viscosity. In
this section and below, constants are independent of $h$. They depend
only on $L$, $\kappa$ and the indicated decay rates.

For a functional $g\in V_1^*$ we write $\Pi_hg$ for the element of $X_h$
that represents the restriction of $g$ to $X_h$ in the weighted pairing.
The discrete dual norm is
\[
 \norm r_{-1,h}=\sup_{0\ne\psi\in X_h}\frac{|\ip r\psi|}{|\psi|_1},
\]
and $\norm{\Pi_hg}_{-1,h}\le\norm g_{V_1^*}$.

\begin{lemma}[Projection, residual and graph recovery]\label{lem:projection}
For $w\in D(A)$ let
$d_h(w)=\Pi_hAw-A_h\Pi_hw$ be the residual. For $w\in V_2\cap H^s$,
$s=2,3,4$,
\begin{equation}\label{eq:approx}
 \norm{w-\Pi_hw}_{H^j}\le Ch^{s-j}\norm w_{H^s},\qquad j=0,1,2,
\end{equation}
and $\Pi_h$ is bounded on $V_1$ and on $V_2$. For
$x\in X_h$ and $w\in D(A)$,
\begin{align}
 \norm x_{H^2}&\le C\norm{A_hx}_H,\label{eq:graph}\\
 \norm{A_h\Pi_hw}_H&\le C\norm w_{H^3}.\label{eq:graphproj}
\end{align}
The residual satisfies
\begin{equation}\label{eq:residual}
 \ip{d_h(w)}\psi=\int_0^L(w-\Pi_hw)''\,\overline{(\rho\psi)'}\dx,
 \qquad
 \begin{cases}
 \norm{d_h(w)}_{-1,h}\le Ch^{s-2}\norm w_{H^s},\\
 \norm{d_h(w)}_H\le Ch^{s-3}\norm w_{H^s},
 \end{cases}\quad s=3,4 .
\end{equation}
\end{lemma}
\begin{proof}
We use Hermite interpolation~\cite[Theorem~3.1.6]{Ciarlet}, the
best-approximation property of $\Pi_h$ and inverse
estimates~\cite[Section~3.2]{Ciarlet} on $\Pi_hw-I_hw=\Pi_h(w-I_hw)$. This
proves
\eqref{eq:approx} and the $V_2$ stability. For the $V_1$ stability we use the cubic interpolant $Q_hw$ that matches nodal values and sets all nodal
slopes to zero. It satisfies $\norm{(Q_hw)'}_2\le C\norm{w'}_2$ and
$\norm{w-Q_hw}_2\le Ch\norm{w'}_2$. The same decomposition applies.

For graph recovery we set $w=A^{-1}A_hx$. Then $\norm w_{H^3}\le
C\norm{A_hx}_H$ and $b(w-x,\psi)=0$ for $\psi\in X_h$. We test with $e=x-I_hw$. Coercivity and interpolation give
$\norm{e'}_2\le C\norm{(w-I_hw)''}_2\le Ch\norm{A_hx}_H$. Now \eqref{eq:graph} follows from an inverse estimate and the $H^2$ stability of $I_h$.
Subtracting the weak equations, we obtain the exact residual identity. The two residual bounds follow from approximation and one inverse estimate. The strong residual bound with $s=3$ gives \eqref{eq:graphproj}.
\end{proof}

\subsection{Exponential stability and Kato smoothing}\label{sec:stability}
Equation~\eqref{eq:energy} holds exactly on $X_h$. As in
Proposition~\ref{prop:rosier}, $B_h$ is invertible and maximal
accretive, and
\begin{equation}\label{eq:disenergy}
 \norm{T_h(t)}_{H\to H}\le e^{-\lambda t},\qquad
 |x|_1^2\le a^{-1}\norm{B_hx}_H\norm x_H,\qquad x\in X_h .
\end{equation}
The same proof yields the mesh-uniform Kato smoothing
\begin{equation}\label{eq:disKato}
 \norm{T_h(t)x}_H^2+2a\int_0^t|T_h(s)x|_1^2\,ds\le\norm x_H^2,
 \qquad x\in X_h,\ t\ge0 .
\end{equation}
The adjoint $B_h^*$ on $X_h$ satisfies
$\Re\ip{B_h^*x}x=\Re\ip{B_hx}x\ge a|x|_1^2$. For $y_t=A_hy+g$ in $X_h$
the argument of \eqref{eq:forced} yields, with the same constant,
\begin{equation}\label{eq:disforced}
 \E_\nu(y)\le C_{L,\nu}
 \big(\norm{y(0)}_H+\norm{e^{\nu t}g}_{L^2_t(-1,h)}\big),
 \qquad 0<\nu<\lambda ,
\end{equation}
where $L^2_t(-1,h)$ is the $L^2$ space in time on $(0,\infty)$ with spatial norm $\norm\cdot_{-1,h}$.

For data of small energy the nonlinear scheme keeps the exponential decay at any rate $\nu<\lambda$, with a radius depending on $\nu$.

\begin{proposition}[Exponential decay for small energy]\label{prop:nonlinstab}
Let $0<\nu<\lambda$ and let $u_0\in H$ be real with
$\norm{u_0}_H\le9\pi(1-\nu/\lambda)/(2L)$. Then for every mesh and every
$t\ge0$,
\[
 \norm{u_h(t)}_H\le e^{-\nu t}\norm{\Pi_hu_0}_H,\qquad
 2\nu p^2\int_0^\infty|u_h|_1^2\,dt\le\norm{\Pi_hu_0}_H^2 .
\]
\end{proposition}
\begin{proof}
For real $w\in H_0^1$,
$\ip{F(w,w)}w=-\int_0^Lw^2w'\rho\dx=\frac13\int_0^Lw^3\dx$, since
$\rho'=1$ and $w$ vanishes at both ends. By the bounds $\norm w_\infty\le L^{1/2}|w|_1$, $\norm w_2\le(L/\pi)|w|_1$ and $\norm w_2\le L^{-1/2}\norm w_H$, we obtain
\[
 \tfrac13\Big|\int_0^Lw^3\dx\Big|\le\frac{L}{3\pi}\norm w_H|w|_1^2 .
\]
Since $\Pi_h$ is the $H$-orthogonal projection,
$\ip{\Pi_hF(u_h,u_h)}{u_h}=\ip{F(u_h,u_h)}{u_h}$. Hence, by \eqref{eq:energy}, we have
\[
 \frac12\frac d{dt}\norm{u_h}_H^2
 \le-\Big(a-\frac L{3\pi}\norm{u_h}_H\Big)|u_h|_1^2 .
\]
Since $a=3/2$, the bracket is at least $a\nu/\lambda=\nu p^2$ as long as
$\norm{u_h}_H\le9\pi(1-\nu/\lambda)/(2L)$. Then the norm does not
increase. The bound on $\norm{u_h}_H$ therefore holds for all $t\ge0$. Integrating, we obtain the second estimate. The first follows from the inequality $\norm{u_h}_H\le p|u_h|_1$ and $\lambda=a/p^2$.
\end{proof}
The same computation applies to the continuous solution whenever the
energy identity is available for it, for instance for data in a small
ball of $V_3$. The radius in Proposition~\ref{prop:nonlinstab} is explicit and proportional to $1/L$, in
accordance with the scaling of \eqref{eq:pde}.

Error equations are linearizations along a pair of solutions. We write
$A_\bullet=A$ in the continuous case and $A_\bullet=A_h$ in the discrete
case. In these equations small solutions enter only as coefficients. We treat the
coefficient terms by the forced estimates \eqref{eq:forced} and
\eqref{eq:disforced} and one absorption. For this we shall repeatedly use
\begin{equation}\label{eq:prodenergy}
 \norm{F(f,g)}_H\le C|f|_1|g|_1,\qquad
 \norm{F(f,g)}_{V_1^*}\le C|f|_1\norm g_H,\qquad f,g\in H_0^1 .
\end{equation}
Both follow from integration by parts and the embedding
$H_0^1\hookrightarrow L^\infty$. The second also holds with $f$ and $g$ interchanged. It also
holds for $F_h$ with the norm $\norm\cdot_{-1,h}$.

\subsection{Hilbert direct sums}\label{sec:directsum}
Uniformity in the mesh will be obtained from statements about a single
generator. Let $S_n$, $n\in\mathbb N$, be contraction semigroups with
generators $G_n$ on Hilbert spaces $\mathcal H_n$, and let
$\mathfrak H=\bigoplus_n\mathcal H_n$ be their $\ell^2$ direct sum.
Let $\iota_n:\mathcal H_n\to\mathfrak H$ be the inclusion as the $n$-th
summand and $\pi_n:\mathfrak H\to\mathcal H_n$ the coordinate map, so that
$\pi_n\iota_n=I$ and $\iota_n\pi_n$ is the orthogonal projection onto the
$n$-th summand.

\begin{lemma}[Direct sums]\label{lem:directsum}
The operators $\mathfrak S(t)=\bigoplus_nS_n(t)$ form a contraction
$C_0$-semigroup on $\mathfrak H$. Its generator is
\begin{equation}\label{eq:diagonal}
 \mathfrak G(x_n)_n=(G_nx_n)_n,\qquad
 D(\mathfrak G)=\Big\{(x_n)\in\mathfrak H:\ x_n\in D(G_n),\
 \sum_n\norm{G_nx_n}^2<\infty\Big\}.
\end{equation}
Moreover, for $z\in\C$ and $t>0$:
\begin{enumerate}[label=\rm(\roman*)]
\item $z\in\rho(\mathfrak G)$ if and only if $z\in\bigcap_n\rho(G_n)$ and
 $\sup_n\norm{R(z,G_n)}<\infty$. In this case
 $\norm{R(z,\mathfrak G)}=\sup_n\norm{R(z,G_n)}$.
\item $\mathfrak S(t)\mathfrak H\subset D(\mathfrak G)$ if and only if
 $S_n(t)\mathcal H_n\subset D(G_n)$ for all $n$ and
 $\sup_n\norm{G_nS_n(t)}<\infty$. In this case
 $\norm{\mathfrak G\mathfrak S(t)}=\sup_n\norm{G_nS_n(t)}$.
\end{enumerate}
\end{lemma}
\begin{proof}
\emph{Step 1: the semigroup.} Each $S_n(t)$ is a contraction. Hence
$\mathfrak S(t)$ is a well-defined contraction. The semigroup law holds
componentwise. For $x=(x_n)\in\mathfrak H$,
\[
 \norm{\mathfrak S(t)x-x}^2=\sum_n\norm{S_n(t)x_n-x_n}^2 .
\]
Each term tends to zero as $t\to0$ and is bounded by $4\norm{x_n}^2$.
Strong continuity follows by dominated convergence. By definition,
$\pi_n\mathfrak S(t)=S_n(t)\pi_n$ and $\mathfrak S(t)\iota_n=\iota_nS_n(t)$.

\emph{Step 2: the generator.} Let $\mathfrak G$ be the generator and let
$\mathfrak G_d$ be the diagonal operator on the set $D_d$ displayed in
\eqref{eq:diagonal}. Fix $z$ with $\Re z>0$. The resolvent is the Laplace
transform $R(z,\mathfrak G)y=\int_0^\infty e^{-zt}\mathfrak S(t)y\,dt$.
Applying the bounded map $\pi_n$ under the integral, we obtain
$\pi_nR(z,\mathfrak G)y=R(z,G_n)\pi_ny$. Let $x\in D(\mathfrak G)$ and
$y=(z-\mathfrak G)x$. Then $x_n=R(z,G_n)y_n\in D(G_n)$ and
$G_nx_n=zx_n-y_n=\pi_n\mathfrak Gx$. Hence
$\sum_n\norm{G_nx_n}^2=\norm{\mathfrak Gx}^2<\infty$. Thus
$D(\mathfrak G)\subset D_d$ and $\mathfrak G=\mathfrak G_d$ on
$D(\mathfrak G)$. Conversely, $z-\mathfrak G$ maps $D(\mathfrak G)$ onto
$\mathfrak H$. On $D_d$ the operator $z-\mathfrak G_d$ is injective, because each
$z-G_n$ is injective. If $x\in D_d$, choose $\tilde x\in D(\mathfrak G)$
with $(z-\mathfrak G)\tilde x=(z-\mathfrak G_d)x$. Then
$(z-\mathfrak G_d)(x-\tilde x)=0$. Hence $x=\tilde x\in D(\mathfrak G)$,
and $\mathfrak G=\mathfrak G_d$.

\emph{Step 3: part (i).} By Step~2, $\iota_n$ maps $D(G_n)$ into
$D(\mathfrak G)$ with $\mathfrak G\iota_n=\iota_nG_n$. Moreover, $\iota_n\pi_n$
maps $D(\mathfrak G)$ into itself and commutes with $\mathfrak G$ there.
Let $z\in\rho(\mathfrak G)$. Then $R(z,\mathfrak G)$ commutes with
$\iota_n\pi_n$. For $y_n\in\mathcal H_n$ the element
$x=R(z,\mathfrak G)\iota_ny_n$ satisfies $x=\iota_n\pi_nx$. Hence
$x=\iota_nx_n$ with $x_n\in D(G_n)$ and $(z-G_n)x_n=y_n$. If
$(z-G_n)x_n=0$, then $(z-\mathfrak G)\iota_nx_n=0$ and $x_n=0$. Hence
$z\in\rho(G_n)$ and $R(z,G_n)=\pi_nR(z,\mathfrak G)\iota_n$, whose norm
is at most $\norm{R(z,\mathfrak G)}$. Conversely, let
$z\in\bigcap_n\rho(G_n)$ with $\norm{R(z,G_n)}\le M$ for all $n$. The
diagonal operator $\mathfrak K=\bigoplus_nR(z,G_n)$ has norm
$\sup_n\norm{R(z,G_n)}\le M$. It maps $\mathfrak H$ into $D(\mathfrak G)$,
because $G_nR(z,G_n)=zR(z,G_n)-I$. It is a two-sided inverse of
$z-\mathfrak G$. The norm identity follows from the two inequalities.

\emph{Step 4: part (ii).} If $\mathfrak S(t)\mathfrak H\subset
D(\mathfrak G)$, then $\mathfrak G\mathfrak S(t)$ is closed and
everywhere defined, hence bounded. Applied to $\iota_nx_n$, it yields
$\iota_nG_nS_n(t)x_n$. Hence $\norm{G_nS_n(t)}\le\norm{\mathfrak G\mathfrak
S(t)}$. Conversely, if the blocks are uniformly bounded, then
$\sum_n\norm{G_nS_n(t)x_n}^2\le\sup_n\norm{G_nS_n(t)}^2\norm x^2$. Thus
$\mathfrak S(t)x\in D(\mathfrak G)$ by \eqref{eq:diagonal}, and the norm
identity follows. The proof is complete.
\end{proof}

\begin{proposition}[Uniform Crandall--Pazy regularity]\label{prop:uniformCP}
Let $\{S_j\}_{j\in\mathfrak J}$ be a family of contraction semigroups with
generators $G_j$ on Hilbert spaces, and let $0<\beta\le1$. The following
are equivalent.
\begin{enumerate}[label=\rm(\alph*)]
\item Every $S_j$ is immediately differentiable and
 $\sup_{j\in\mathfrak J}\sup_{0<t\le1}t^{1/\beta}\norm{G_jS_j(t)}<\infty$.
\item There are $\Omega,C>0$ such that $i\omega\in\rho(G_j)$ and
 $\norm{R(i\omega,G_j)}\le C|\omega|^{-\beta}$ for all $j\in\mathfrak J$
 and all $|\omega|\ge\Omega$.
\end{enumerate}
\end{proposition}
\begin{proof}
(b)$\Rightarrow$(a). By Theorem~\ref{thm:wakaiki} every $S_j$ is in
$\CP_\beta$. Suppose that the supremum in (a) is infinite. Then there are
$j_n\in\mathfrak J$ and $t_n\in(0,1]$ with
$t_n^{1/\beta}\norm{G_{j_n}S_{j_n}(t_n)}\ge n$. Let $\mathfrak S$ be the
direct sum of $S_{j_1},S_{j_2},\dots$. By (b) and
Lemma~\ref{lem:directsum}(i),
$\norm{R(i\omega,\mathfrak G)}\le C|\omega|^{-\beta}$ for
$|\omega|\ge\Omega$. By Theorem~\ref{thm:wakaiki}, $\mathfrak S\in\CP_\beta$: there are $t_*\in(0,1]$ and $C'$ with
$\norm{\mathfrak G\mathfrak S(t)}\le C't^{-1/\beta}$ for $0<t\le t_*$. For
$t_*\le t\le1$, $\mathfrak G\mathfrak S(t)=\mathfrak S(t-t_*)\mathfrak
G\mathfrak S(t_*)$ yields
$\norm{\mathfrak G\mathfrak S(t)}\le C't_*^{-1/\beta}\le
C't_*^{-1/\beta}t^{-1/\beta}$. By Lemma~\ref{lem:directsum}(ii), the
same bound holds for every $G_{j_n}S_{j_n}(t)$, a contradiction for large
$n$.

(a)$\Rightarrow$(b). Suppose (b) fails. Then there are $j_n$ and
$\omega_n$ with $|\omega_n|\ge n$ such that either
$i\omega_n\notin\rho(G_{j_n})$ or
$|\omega_n|^\beta\norm{R(i\omega_n,G_{j_n})}\ge n$. By (a) and
Lemma~\ref{lem:directsum}(ii), the direct sum of the $S_{j_n}$ belongs to
$\CP_\beta$. By Theorem~\ref{thm:wakaiki} there are $\Omega',C'$ with
$\norm{R(i\omega,\mathfrak G)}\le C'|\omega|^{-\beta}$ for
$|\omega|\ge\Omega'$. Lemma~\ref{lem:directsum}(i) transfers this bound to
every $G_{j_n}$, which contradicts the choice of $(j_n,\omega_n)$ for
$n>\max\{\Omega',C'\}$.
\end{proof}
Only countable direct sums are used. The index set may be uncountable,
for instance the set of all admissible meshes. The analogous
characterization of uniform exponential stability is due to Liu and
Zheng~\cite{LiuZheng94}. It is a uniform version of the
Gearhart--Huang--Pr\"uss theorem~\cite{Gearhart,Huang,Pruss}.

\subsection{Uniform resolvent decay and Crandall--Pazy regularity}
\label{sec:uniformcp}
In this subsection, we prove the mesh-uniform resolvent bound and the
Crandall--Pazy regularity of the W-scheme. The main result reads as
follows.

\begin{theorem}[Uniform resolvent decay and Crandall--Pazy regularity]\label{thm:cp}
For all meshes of Section~\ref{sec:scheme}, with constants independent of
$h$,
\begin{align}
 \norm{(i\tau-A_h)^{-1}}_{H\to H}&\le C\langle\tau\rangle^{-2/3},
 \qquad\tau\in\R,\label{eq:resolvent}\\
 \norm{B_h^\alpha T_h(t)}_{H\to H}&\le Ct_\wedge^{-3\alpha/2}e^{-\lambda t},
 \qquad 0\le\alpha\le1,\ t>0,\label{eq:cp}\\
 |T_h(t)x|_1&\le Ct_\wedge^{-3/4}e^{-\lambda t}\norm x_H,
 \qquad x\in X_h,\label{eq:cpH1}\\
 \norm{T_h(t)\Pi_hg}_H&\le Ct_\wedge^{-3/4}e^{-\lambda t}\norm g_{V_1^*},
 \qquad g\in V_1^*.\label{eq:weaksmoothing}
\end{align}
The continuous operator satisfies the same estimates, with $\Pi_h$
omitted. The exponent $3/2$ in \eqref{eq:cp} with $\alpha=1$ cannot be
improved uniformly in $h$.
\end{theorem}
\begin{proof}
\emph{Resolvent.} Let $x\in X_h$, $\norm x_H=1$, $f=(i\tau-A_h)x$,
$r_f=\norm f_H$ and $k=|x|_1$. Equation~\eqref{eq:energy} holds on
$X_h$. Taking the real part of $\ip fx$, we obtain $ak^2\le r_f$.
Together with $1\le pk$, this yields $r_f\ge\lambda$. Fix a small
$\gamma\in(0,L/(8p)]$ and put $\ell=\gamma/k\le L/8$. If $h\ge\ell$, then, by inverse estimates, $\norm{x''}_2\le Ch^{-1}k\le C\gamma^{-1}k^2$. Hence
$|b(x,x)|\le C_\gamma k^3$. Testing with $x$ then yields
$|\tau|\le r_f+C_\gamma k^3\le Cr_f^{3/2}$.

Suppose $h<\ell$. Let $z\in C_c^\infty(0,L)$ be the function of
Section~\ref{sec:cp}, which satisfies \eqref{eq:mollify}, and take
$\psi=I_hz\in X_h$. On each element $(I_hz)''$ is the $L^2$
projection of $z''$ onto linear polynomials. Indeed, the interpolation error and
its first derivative vanish at both element ends. Hence its second
derivative is orthogonal to all linear functions. Therefore
\[
 \norm{x-\psi}_H\le C(\ell k+h^2\ell^{-1}k)\le2C\gamma,
 \qquad \norm{\psi''}_2\le C\ell^{-1}k,
 \qquad \norm{\psi'}_2\le Ck .
\]
Choose $\gamma$ so that the first bound is at most $1/2$. Then
$|\ip x\psi|\ge1/2$. Since $\psi(0)=\psi'(0)=0$ and $x'(L)=0$, one more
integration by parts yields
\[
 b(x,\psi)=\int_0^Lx'\,\overline{(\rho\psi)''}\dx,
 \qquad |b(x,\psi)|\le C_\gamma k^3 .
\]
Testing $f$ against $\psi$, we obtain
$\frac12|\tau|\le r_f\norm\psi_H+|b(x,\psi)|\le C(r_f+r_f^{3/2})
\le Cr_f^{3/2}$, where the last step uses $r_f\ge\lambda$. Hence
$r_f\ge c\max\{\lambda,|\tau|^{2/3}\}$, which proves \eqref{eq:resolvent}.

\emph{Crandall--Pazy regularity.} By \eqref{eq:resolvent} the family
$\{T_h\}$ satisfies condition (b) of Proposition~\ref{prop:uniformCP}
with $\beta=2/3$. Hence
$\sup_h\sup_{0<t\le1}t^{3/2}\norm{B_hT_h(t)}<\infty$. For $t\ge1$ we write $B_hT_h(t)=T_h(t-1)B_hT_h(1)$ and use \eqref{eq:disenergy}. This proves
\eqref{eq:cp} for $\alpha=1$. The moment inequality for maximal accretive
operators has constants independent of $h$, which yields the case $0<\alpha<1$.
Applying \eqref{eq:disenergy} to $T_h(t)x$, we obtain \eqref{eq:cpH1}.

\emph{Weak sources.} On $X_h$, $B_h$ commutes with $T_h$, so
$\norm{B_h^*T_h(t)^*}=\norm{B_hT_h(t)}$. Combining this with $a|y|_1^2\le\Re\ip{B_h^*y}y$, \eqref{eq:cp} and the contraction property, we obtain
$|T_h(t)^*\psi|_1\le Ct_\wedge^{-3/4}e^{-\lambda t}\norm\psi_H$ for
$\psi\in X_h$. Since
$\ip{T_h(t)\Pi_hg}\psi=\langle g,T_h(t)^*\psi\rangle$, \eqref{eq:weaksmoothing} follows by duality. For the continuous operator we use the adjoint
semigroup and \eqref{eq:adjoint}. The last assertion of the theorem is proved below.
\end{proof}
In particular the family $\{T_h\}$ is uniformly of class $\CP_{2/3}$.
We now prove the last assertion of Theorem~\ref{thm:cp}. By
Theorem~\ref{thm:ctcp}(iii) the index $3/2$ is already optimal for the
continuous semigroup. Suppose that $\norm{B_hT_h(t)}\le Ct^{-\gamma}$ for $0<t\le1$ and
all meshes, for some $\gamma<3/2$. Fix $f\in H$ and $0<t\le1$, and put
$w_h=B_hT_h(t)\Pi_hf$. Then $B_h^{-1}\Pi_hw_h=T_h(t)\Pi_hf$, which tends
to $T(t)f$ as $h\to0$ by Proposition~\ref{thm:positive}.
By Lemma~\ref{lem:resdefect} at $z=0$ and \eqref{eq:approx}, we have
$\norm{B_h^{-1}\Pi_h-B^{-1}}\le Ch$. The family $(w_h)$ is bounded in $H$. Therefore every weak
limit point $w$ of $w_h$ satisfies $B^{-1}w=T(t)f$, that is,
$w=BT(t)f$. Hence $\norm{AT(t)}\le Ct^{-\gamma}$ for $0<t\le1$, which
contradicts Theorem~\ref{thm:ctcp}(iii). The last assertion is not used
elsewhere. Its proof uses Proposition~\ref{thm:positive} and
Lemma~\ref{lem:resdefect}, which do not depend on it.

\subsection{Analyticity of the discrete solution map}\label{sec:uniformanalytic}
In this subsection, we show that the analyticity of
Theorem~\ref{thm:continuous} carries over to the discrete
solution maps, with constants independent of $h$.

\begin{proposition}[Mesh-uniform analyticity]\label{prop:uniformanalytic}
Fix $0<\nu<\lambda$. There is $r>0$, depending only on $L$ and $\nu$,
such that for every mesh the discrete solution map $f\mapsto u_h$,
$u_h(0)=f\in X_h$, extends to a holomorphic map from the complex ball
$\{f\in X_h:\norm f_H<r\}$ into $\X_0$. There are homogeneous
polynomials $P_{n,h}$ of degree $n$ with
\[
 u_h=\sum_{n\ge1}P_{n,h}(f),\qquad P_{1,h}(f)=T_h(\cdot)f,\qquad
 \norm{P_{n,h}(f)}_{\X_0}\le4^{n-1}q^{n-1}l^n\norm f_H^n ,
\]
where $l$ and $q$ are the bounds of $\LL$ and $\QQ$ in
\eqref{eq:lowbilinear}. The series converges absolutely in $\X_0$, and
the Taylor remainder estimate \eqref{eq:taylor} holds for
$u_h$ with constants independent of $h$.
\end{proposition}
\begin{proof}
Equation~\eqref{eq:energy} holds exactly on $X_h$. The weak product bound
and Agmon's inequality use only $H_0^1$, and
$\norm{\Pi_hg}_{-1,h}\le\norm g_{V_1^*}$. Hence \eqref{eq:disforced} and
the bilinear bound \eqref{eq:lowbilinear} hold for $T_h$ and $F_h$ with
the same constants. The proof of Theorem~\ref{thm:continuous} then applies
verbatim. Small discrete solutions are unique, so the sum of the series is
$u_h$.
\end{proof}

\subsection{Resolvent comparison and a common contour}
\label{sec:spectral}
In this subsection, we compare the discrete resolvents with the continuous
one. Lemma~\ref{lem:CPcontour}
provides a contour representation of the semigroups with common bounds,
which is used in Section~\ref{sec:samenorm}.

\begin{lemma}[Uniform resolvent defect]\label{lem:resdefect}
Let $C_R$ be a constant such that, by Theorem~\ref{thm:cp}, $\norm{R(i\omega,A)}$ and $\norm{R(i\omega,A_h)}$ are at most
$C_R\langle\omega\rangle^{-2/3}$, and let
$c_\Gamma=\min\{\frac14,(2C_R)^{-1}\}$. For $z=-\sigma+i\omega$ with
$0\le\sigma\le c_\Gamma\langle\omega\rangle^{2/3}$ and $0\le\theta\le3$,
\[
 \norm{\{R(z,A_h)\Pi_h-\Pi_hR(z,A)\}B^{-\theta/3}}_{H\to H}
 \le Ch\langle\omega\rangle^{-\theta/3}.
\]
\end{lemma}
\begin{proof}
Both resolvents exist with norm at most $2C_R\langle\omega\rangle^{-2/3}$,
by a Neumann series around $i\omega$, as in \eqref{eq:strip}. Let
$w=R(z,A)f$. Since $(z-A_h)\Pi_hw=\Pi_hf+d_h(w)$,
\[
 y:=R(z,A_h)\Pi_hf-\Pi_hR(z,A)f=-R(z,A_h)d_h(w).
\]
For $\psi\in X_h$ put $x=R(\bar z,A_h^*)\psi$, where $A_h^*$ is the
$H$-adjoint of $A_h$ on $X_h$. Then $\ip y\psi=-\ip{d_h(w)}x$, so
$|\ip y\psi|\le\norm{d_h(w)}_{-1,h}|x|_1$. Taking the real part of $\ip{(\bar z-A_h^*)x}x=\ip\psi x$ and using \eqref{eq:energy}, we obtain
\[
 a|x|_1^2\le\norm\psi_H\norm x_H+\sigma\norm x_H^2
 \le C\langle\omega\rangle^{-2/3}\norm\psi_H^2 .
\]
By \eqref{eq:residual} with $s=3$,
$\norm{d_h(w)}_{-1,h}\le Ch\norm w_{H^3}\le Ch\norm{AR(z,A)f}_H$.
Finally $\norm{AR(z,A)}\le C\langle\omega\rangle^{1/3}$ and
$\norm{AR(z,A)B^{-1}}=\norm{R(z,A)}\le C\langle\omega\rangle^{-2/3}$.
By complex interpolation on the scale $V_\theta$, we obtain
$\norm{AR(z,A)B^{-\theta/3}}\le C\langle\omega\rangle^{1/3-\theta/3}$.
The three bounds together give the claim.
\end{proof}

\begin{lemma}[A common contour with absolutely convergent integrals]\label{lem:CPcontour}
There is $c_\Gamma'\in(0,c_\Gamma]$, depending only on $L,\kappa$, such
that the upward-oriented curve
$\Gamma:\ z(\omega)=-c_\Gamma'\langle\omega\rangle^{2/3}+i\omega$
lies in the resolvent sets of $G=A$ and of every $G=A_h$, and
\[
 \norm{R(z(\omega),G)}\le C\langle\omega\rangle^{-2/3},\qquad
 e^{tG}=\frac1{2\pi i}\int_\Gamma e^{zt}R(z,G)\,dz,\quad t>0 .
\]
The integral converges in operator norm. For $A_h$ the formula holds on
$X_h$ and may then be composed with $\Pi_h$.
\end{lemma}
\begin{proof}
The strip bound of Lemma~\ref{lem:resdefect} covers $\Gamma$ at large
$|\omega|$. At bounded frequencies we choose the leftward displacement
smaller than $\lambda/2$ and use the contraction property. Thus the region between
$\Gamma$ and any vertical line $\Re z=b>0$ is free of spectrum, and the
resolvent decays there as $|\Im z|\to\infty$. We decrease $c_\Gamma'$ further so
that $|z'(\omega)|\le2$. The integral defines a norm-continuous
family of operators $S_G(t)$, $t>0$, and
\[
 \norm{S_G(t)}\le C\int_{\R}e^{-c_\Gamma't\langle\omega\rangle^{2/3}}
 \langle\omega\rangle^{-2/3}\,d\omega\le C(1+t^{-1/2})e^{-c't}.
\]
In particular its Laplace transform exists for real $s_0>0$. Fubini's
theorem is justified by $\norm{R(z,G)}/(s_0-\Re z)=O(|\omega|^{-4/3})$.
After the time integration the integrand is $R(z,G)/(s_0-z)$, of order
$|\omega|^{-5/3}$ on $\Gamma$.

To identify $S_G(t)$, we multiply by $G^{-1}$, whose norm is at most
$\lambda^{-1}$, and deform the integral of $R(z,G)G^{-1}/(s_0-z)$ to a
vertical line $\Re z=b>s_0$. The integrals over the horizontal joining segments at heights
$\pm R_0$ tend to zero as $R_0\to\infty$. Indeed, the segments have length $O(R_0^{2/3})$, while the
integrand is $O(R_0^{-5/3})$. The only enclosed pole is at $z=s_0$.
Finally we use $R(z,G)G^{-1}=z^{-1}(R(z,G)+G^{-1})$ and
$\norm{R(b+i\omega,G)}\le(b+\lambda)^{-1}$. The regularized integral on
$\Re z=b$ is $O(b^{-1})$ as $b\to\infty$, since
$\int_{\R}|b+i\omega|^{-1}|b-s_0+i\omega|^{-1}d\omega=O(b^{-1})$.
By Cauchy's formula, we now obtain
$\int_0^\infty e^{-s_0t}S_G(t)G^{-1}dt=R(s_0,G)G^{-1}$. By uniqueness of the
Laplace transform, $S_G(t)G^{-1}=e^{tG}G^{-1}$ for $t>0$. The
range of $G^{-1}$ is dense, so the formula follows. All the estimates
just used have common constants. The proof is complete.
\end{proof}

\section{Error estimates}\label{sec:samenorm}
\setcounter{equation}{0}

In this section, we recall the ordering error and establish the error
estimates.

\subsection{Two orders of approximation}\label{sec:orders}
Let $u$ solve \eqref{eq:pde}, let $u_h$ solve the W-scheme
\eqref{eq:Wscheme}, and let
\[
 v=T(\cdot)u_0,\qquad v_h=T_h(\cdot)\Pi_hu_0
\]
be the linear solutions. The linear discrete simulation $v_h$ can be
reached in two orders. Each order divides the total error in its own
way:
\begin{equation}\label{eq:allocation}
 u-v_h=\underbrace{(u-v)}_{\text{linearization}}
      +\underbrace{(v-v_h)}_{\text{discretization}}
      =\underbrace{(u-u_h)}_{\text{discretization}}
      +\underbrace{(u_h-v_h)}_{\text{linearization}} .
\end{equation}
Changing the order changes the linearization share from $u-v$ to
$u_h-v_h$, and the discretization share from $v-v_h$ to $u-u_h$. The two changes cancel in the total error. As in the introduction, the
ordering error is
\begin{equation}\label{eq:defect}
 \Cd_h=(u-v)-(u_h-v_h)=(u-u_h)-(v-v_h).
\end{equation}
It measures how much the nonlinear correction
predicted by the discrete model differs from the continuous one.
Equivalently, it measures how much the nonlinearity changes the discretization
error. We call the discretization shares $v-v_h$ and $u-u_h$ the
\emph{linear error} and the \emph{full error}.

Put $e=u-u_h$ and $m=\frac12(u+u_h)$, and write
\[
 L_h(t)=T_h(t)\Pi_h-T(t),\quad (L_h*g)(t)=\int_0^tL_h(t-s)g(s)\,ds,
 \quad \QQ_h(f,g)(t)=\int_0^tT_h(t-s)F_h(f(s),g(s))\,ds .
\]
Then $v_h-v=L_h(\cdot)u_0$ and, since $\N(u)-\N(u_h)=2F(e,m)$,
\begin{align}
 e&=-L_h(\cdot)u_0+\Cd_h,\qquad e(0)=(I-\Pi_h)u_0,\qquad \Cd_h(0)=0,
 \label{eq:splitting}\\
 \Cd_h&=-L_h*\N(u)+2\QQ_h(e,m).\label{eq:Cdmild}
\end{align}
The linear term $L_h(\cdot)u_0$ is absent from $\Cd_h$. This cancellation
is behind every estimate of the ordering error below. Note that
$L_h(t)u_0=-T(t)(I-\Pi_h)u_0+(T_h(t)-T(t))\Pi_hu_0$ contains the initial
projection error propagated by $T$ and a consistency error.

Throughout, $\varepsilon$ denotes a norm of the datum, and no direction
is fixed. If a direction $\phi$ is fixed and $u_0=\delta\phi$, the ordering error can be read as a
mixed difference. Let $f(\delta,h)=u_h(\delta\phi)$
and $f(\delta,0)=u(\delta\phi)$, put $q(\delta,h)=f(\delta,h)/\delta$
for $\delta\ne0$, and extend $q$ to $\delta=0$ by the derivative of $f$
with respect to $\delta$ at $\delta=0$. Then
\[
 \Cd_h=\delta\{q(\delta,0)-q(0,0)-q(\delta,h)+q(0,h)\}.
\]
Thus $\Cd_h/\delta$ is a mixed second difference of the normalized
solution. For $t\ge t_0>0$, bounds for the four shares give only
$\Cd_h=O(\min\{\delta^2,h^p|\delta|\})$, with constants depending on
$t_0$. The product bound
$O(h^p\delta^2)$ requires the cancellation in \eqref{eq:Cdmild}.

\subsection{Statement of the error estimates}\label{sec:mainsame}
In Theorems~\ref{thm:linear} and~\ref{thm:ordering} the error and the
datum are measured in the same norm, the norm of $H$. The constants
depend only on $L$, $\kappa$ and the decay rate. We write
$\varepsilon=\norm{u_0}_H$. For sufficiently small $\varepsilon$ the four
solutions exist for all $t\ge0$, by Theorem~\ref{thm:continuous} and
Proposition~\ref{prop:nonlinstab}. Theorem~\ref{thm:DA} shows what a
higher norm of the datum adds.

Pointwise rates for the full error cannot hold uniformly on a ball down to
$t=0$. The kernel of $\Pi_h$ contains nonzero smooth compactly supported
functions. Hence $\norm{e(0)}_H=\norm{(I-\Pi_h)u_0}_H$ can
equal the radius of the ball for every $h$.

\begin{theorem}[Linear and full errors, linearization shares]\label{thm:linear}
Fix $L,\kappa$ and $0<\nu<\lambda$. There are $r,C>0$, depending only
on $L,\kappa,\nu$, with the following property. Let $u_0\in H$ be real
with $\varepsilon=\norm{u_0}_H\le r$, and let $0<h\le\min\{1,L\}$. Then,
for $t>0$,
\begin{align}
 \norm{(v-v_h)(t)}_H
 &\le Ce^{-\nu t}\min\{1,ht_\wedge^{-3/2}\}\,\varepsilon,\label{eq:sn-linear}\\
 \norm{(u-u_h)(t)}_H
 &\le Ce^{-\nu t}\big\{\min\{1,ht_\wedge^{-3/2}\}\,\varepsilon
 +h^{1/6}\varepsilon^2\big\}.\label{eq:sn-full}
\end{align}
In particular, for every $t_0>0$, the errors converge to zero uniformly
for $t\ge t_0$ and $\varepsilon\le r$. The two linearization shares
satisfy
\begin{equation}\label{eq:sn-horizontal}
 \sup_{t\ge0}e^{\nu t}\norm{(u-v)(t)}_H\le C\varepsilon^2,
 \qquad
 \sup_{t\ge0}e^{\nu t}\norm{(u_h-v_h)(t)}_H\le C\varepsilon^2 .
\end{equation}
\end{theorem}

\begin{theorem}[The ordering error on an energy ball]\label{thm:ordering}
Fix $L,\kappa$ and $0<\nu<\lambda$. There are $r,C>0$, depending only
on $L,\kappa,\nu$, with the following property. Let $u_0\in H$ be real
with $\varepsilon=\norm{u_0}_H\le r$, and let $0<h\le\min\{1,L\}$. The
ordering error satisfies $\Cd_h(0)=0$ and
\begin{equation}\label{eq:sn-ordering}
 \sup_{t\ge0}e^{\nu t}\norm{\Cd_h(t)}_H\le Ch^{1/6}\varepsilon^2 .
\end{equation}
Moreover, let $D_h(u_0)=P_2(u_0)-P_{2,h}(\Pi_hu_0)$ be the discretization
error of the quadratic Taylor coefficient. Then
\[
 \sup_{t\ge0}e^{\nu t}\norm{\Cd_h(t)-D_h(u_0)(t)}_H
 \le Ch^{1/6}\varepsilon^3,\qquad
 \norm{D_h(u_0)(t)}_H\le Ch^{1/2}t_\wedge^{-5/2}e^{-\nu t}
 \varepsilon^2,\quad t>0 .
\]
\end{theorem}

The linearization shares are of order $\varepsilon^2$. Their
difference, the ordering error, carries an additional
mesh factor. In the energy norm this holds uniformly on the whole
half-line, because the linear error $v-v_h$, which contains the initial
projection error, drops out of the mild formula for $\Cd_h$.

A higher norm of the datum buys uniformity down to $t=0$. The next
theorem is the only main result in which the right side carries a higher
norm than the left.

\begin{theorem}[Data in $D(A)$]\label{thm:DA}
Fix $L,\kappa$ and $0<\nu<\lambda$. There are $r,C>0$, depending only
on $L,\kappa,\nu$, such that the following hold for real data and
$0<h\le\min\{1,L\}$. They are independent of $h$, $\hat h$, $u_0$ and
$t$.
\begin{enumerate}[label=\rm(\roman*)]
\item If $u_0\in D(A)$, then $\E_\nu(v-v_h)\le Ch^2\norm{Au_0}_H$.
\item If $u_0\in D(A)$ and $\norm{u_0}_H\le r$, then
\[
 \E_\nu(v-v_h)+\E_\nu(u-u_h)\le Ch^2\norm{u_0}_{H^3}.
\]
\item If $u_0\in D(A)$ and $\varepsilon=\norm{u_0}_H\le r$, then
\[
 \E_\nu(\Cd_h)\le Ch^2\norm{Au_0}_H^{3/2}\varepsilon^{1/2}.
\]
In particular, for all $t\ge0$,
\[
 \norm{\Cd_h(t)}_H\le Ce^{-\nu t}h^2\norm{Au_0}_H^{3/2}\varepsilon^{1/2}.
\]
Moreover, let $h\le\hat h$ be the sizes of two such meshes, and let
$\tilde u=v_h+(u_{\hat h}-v_{\hat h})$. Then
$u-\tilde u=(v-v_h)+\Cd_{\hat h}$, and for all $t\ge0$,
\[
 e^{\nu t}\norm{(u-\tilde u)(t)}_H
 \le Ch^2\norm{Au_0}_H+C\hat h^2\norm{Au_0}_H^{3/2}\varepsilon^{1/2}.
\]
\end{enumerate}
\end{theorem}
In (ii) the constant does not depend on $\norm{u_0}_{H^3}$, which enters
only linearly. Smallness is required only in $H$.

The rates are not claimed to be optimal. Theorems~\ref{thm:linear}
and~\ref{thm:ordering} are proved in Section~\ref{sec:energyball}, and
Theorem~\ref{thm:DA} in Section~\ref{sec:convproof}.

\subsection{Linear errors}\label{sec:linearerror}
The theorems are proved from linear error bounds, some of which carry a
higher norm of the datum. We collect them here.

\begin{proposition}[The linear error at positive times]\label{thm:positive}
Let $0\le\theta\le3$, $f\in V_\theta$, $0<\nu<\lambda$ and $t>0$. Then
\[
 \norm{T_h(t)\Pi_hf-T(t)f}_H^2
 =\norm{(I-\Pi_h)T(t)f}_H^2+\norm{T_h(t)\Pi_hf-\Pi_hT(t)f}_H^2,
\]
with
\begin{align*}
 \norm{(I-\Pi_h)T(t)f}_H
 &\le Ch^3t_\wedge^{-(3-\theta)/2}e^{-\lambda t}\norm f_{V_\theta},\\
 \norm{T_h(t)\Pi_hf-\Pi_hT(t)f}_H
 &\le Ch\,\varsigma_\theta(t)\,e^{-\nu t}\norm f_{V_\theta},
 \qquad
 \varsigma_\theta(t)=\begin{cases}t_\wedge^{-(3-\theta)/2},&\theta<3,\\
                  \log(e/t_\wedge),&\theta=3.\end{cases}
\end{align*}
Here $C$ depends also on $\theta$. In particular
$\norm{L_h(t)}_{H\to H}\le Cht_\wedge^{-3/2}e^{-\nu t}$.
\end{proposition}
\begin{proof}
The two components are orthogonal. For the first, by \eqref{eq:approx}, we have $\norm{(I-\Pi_h)w}_H\le Ch^3\norm{Aw}_H$ on $D(A)$, and
$\norm{AT(t)B^{-\theta/3}}\le Ct_\wedge^{-(3-\theta)/2}e^{-\lambda t}$ by
Theorem~\ref{thm:ctcp}(ii).

For the second, we use the common contour $\Gamma$ of
Lemma~\ref{lem:CPcontour}. Subtracting the two absolutely convergent representations and using Lemma~\ref{lem:resdefect} with $|z'(\omega)|\le2$, we obtain
\[
 \norm{\{T_h(t)\Pi_h-\Pi_hT(t)\}B^{-\theta/3}}
 \le Ch\int_{\R}e^{-c_\Gamma't\langle\omega\rangle^{2/3}}
 \langle\omega\rangle^{-\theta/3}\,d\omega .
\]
For $0<t\le1$ the substitution $s=t\omega^{2/3}$ shows that the integral is at most
$C\varsigma_\theta(t)$. For $t>1$ let $n=\lceil t\rceil$, $t_1=t/n$, and
telescope
\[
 T_h(t)\Pi_h-\Pi_hT(t)=\sum_{j=0}^{n-1}T_h((n-1-j)t_1)
 \{T_h(t_1)\Pi_h-\Pi_hT(t_1)\}T(jt_1).
\]
Both semigroups decay at the rate $\lambda$, and $T$ commutes with
$B^{\theta/3}$. This yields the factor $ne^{-\lambda(n-1)t_1}\le
Ce^{-\nu t}$.
\end{proof}
The weight $t_\wedge^{-3/2}$ is the time cost of the three derivatives
used in the residual bound \eqref{eq:residual}.

For data in $D(A)$ the residual can be estimated by means of time-integrated
$H^4$ regularity. This requires only three derivatives of the datum.

\begin{lemma}[Energy-trajectory errors for data in $D(A)$]
\label{lem:DAenergy}
Let $0<\nu<\lambda$ and $g\in D(A)$. Then
\begin{equation}\label{eq:H4integrated}
 \sup_{t\ge0}e^{\nu t}\norm{T(t)g}_{H^3}
 +\norm{e^{\nu t}T(\cdot)g}_{L^2_tH^4}\le C\norm{Ag}_H,
 \qquad
 \E_\nu(L_h(\cdot)g)\le Ch^2\norm{Ag}_H .
\end{equation}
\end{lemma}
\begin{proof}
Put $v=T(\cdot)g$. Then $Av=T(\cdot)Ag$, and it follows from \eqref{eq:forced} with zero source, applied to $Av$, that
$\int_0^\infty e^{2\nu t}|Av|_1^2\,dt\le C\norm{Ag}_H^2$. Since
$v'''=-Av$ and $v''''=-(Av)'$, this is the integrated $H^4$ bound.
The $H^3$ bound follows from the contraction property. Write
$-L_h(\cdot)g=v-v_h=\eta_v+\Xi$, with $v_h=T_h(\cdot)\Pi_hg$,
$\eta_v=v-\Pi_hv$ and $\Xi=\Pi_hv-v_h$. Then
$\Xi_t=A_h\Xi+d_h(v)$ and $\Xi(0)=0$. By \eqref{eq:disforced} and
the weak residual bound \eqref{eq:residual} with $s=4$,
$\E_\nu(\Xi)\le Ch^2\norm{e^{\nu t}v}_{L^2_tH^4}\le Ch^2\norm{Ag}_H$.
By \eqref{eq:approx},
$\E_\nu(\eta_v)\le Ch^3(\sup_te^{\nu t}\norm v_{H^3}
+\norm{e^{\nu t}v}_{L^2_tH^4})\le Ch^3\norm{Ag}_H$.
\end{proof}

\begin{lemma}[The time kernel of the linear error]\label{lem:linear}
Fix $0<\nu<\lambda$ and put $\Lambda_h(t)=e^{\nu t}\norm{L_h(t)}_{H\to H}$.
Then
\begin{equation}\label{eq:kernels}
 \norm{\Lambda_h}_{L^{p_0}(0,\infty)}\le Ch^{2/(3p_0)},\qquad 1\le p_0\le4 .
\end{equation}
\end{lemma}
\begin{proof}
Fix $\nu<\nu'<\lambda$. Then
$\Lambda_h(t)\le Ce^{-(\nu'-\nu)t}\min\{1,ht_\wedge^{-3/2}\}$, by the
contraction property and Proposition~\ref{thm:positive}. Splitting the integral at $t=h^{2/3}$ and at $t=1$, we obtain
\[
 \int_0^\infty\Lambda_h\,dt
 \le C\Big(h^{2/3}+h\int_{h^{2/3}}^1t^{-3/2}\,dt+h\Big)
 \le Ch^{2/3}.
\]
Since $\Lambda_h\le2$ by the contraction property, $\int_0^\infty\Lambda_h^{p_0}\,dt\le
2^{p_0-1}\int_0^\infty\Lambda_h\,dt$, which gives \eqref{eq:kernels}.
\end{proof}

\subsection{The energy ball}\label{sec:energyball}
Throughout this subsection, $0<\nu<\lambda$ is fixed. For $\phi\in H$ put
$\Phi=T(\cdot)\phi$, $\Phi_h=T_h(\cdot)\Pi_h\phi$ and
\begin{equation}\label{eq:Dh}
 D_h(\phi)=\QQ(\Phi,\Phi)-\QQ_h(\Phi_h,\Phi_h)
 =P_2(\phi)-P_{2,h}(\Pi_h\phi),
\end{equation}
the discretization error of the quadratic Taylor coefficient. By
\eqref{eq:lowbilinear} and its discrete analogue, $D_h(\phi)$ is defined
for every $\phi\in H$ and is bounded uniformly in $h$. Its inputs are
linear evolutions. Hence the smoothing of Section~\ref{sec:discrete}
applies to them directly.

\begin{proposition}[The quadratic coefficient]\label{prop:energyquadratic}
For all real or complex $\phi\in H$ and $t>0$,
\begin{equation}\label{eq:energyquadratic}
 \norm{D_h(\phi)(t)}_H\le Ch^{1/2}t_\wedge^{-5/2}e^{-\nu t}\norm\phi_H^2 .
\end{equation}
\end{proposition}
\begin{proof}
Put $R=\norm\phi_H$, $\nu_1=(\lambda+\nu)/2$ and $\Delta_h=\Phi-\Phi_h$.
By Proposition~\ref{thm:positive}, at the indices $0$ and $2$, we have
\[
 \norm{L_h(r')}_{H\to H}\le Chr_\wedge'^{-3/2}e^{-\nu_1r'},\qquad
 \norm{L_h(r')}_{V_2\to H}\le Chr_\wedge'^{-1/2}e^{-\nu_1r'} .
\]
By Theorem~\ref{thm:ctcp}(ii), $\norm{\Phi(s)}_{V_3}\le Cs_\wedge^{-3/2}e^{-\lambda s}R$. The common
energy identity, with $a-\nu_1p^2>0$, yields
\[
 \int_0^\infty e^{2\nu_1s}(|\Phi(s)|_1+|\Phi_h(s)|_1)^2\,ds\le CR^2 .
\]
By symmetry and bilinearity of $F$, we obtain the exact decomposition
\[
 D_h=-\int_0^tL_h(t-s)F(\Phi,\Phi)\,ds
     +\int_0^tT_h(t-s)\Pi_hF(\Delta_h,\Phi+\Phi_h)\,ds=:I_1+I_2 .
\]
We split each integral at $t/2$.

\emph{Early part of $I_1$.} By Agmon's inequality, $\norm{F(\Phi,\Phi)}_H\le C\norm\Phi_H^{1/2}|\Phi|_1^{3/2}$. By H\"older's
inequality in time,
\[
 \int_0^\tau\norm{F(\Phi,\Phi)}_H\,ds\le CR^{1/2}\tau^{1/4}
 \Big(\int_0^\tau|\Phi|_1^2\,ds\Big)^{3/4}\le C\tau^{1/4}R^2 ,
 \qquad 0<\tau\le1 .
\]
The time-integrated gradient enters with the power $3/4$. Hence this bound is uniform on
the ball. For $t\le1$ the early part is at most $Cht^{-5/4}R^2$. For $t\ge1$ the same argument yields
$\int_0^\infty e^{\nu s}\norm{F(\Phi,\Phi)}_H\,ds\le CR^2$. After
extracting $e^{-\nu t}$ from the convolution, the early part is bounded
by $Cht_\wedge^{-5/4}e^{-\nu t}R^2$ for all $t$.

\emph{Late part of $I_1$.} By \eqref{eq:product} at $\theta=3$,
$\norm{F(\Phi,\Phi)}_{V_2}\le Cs_\wedge^{-3}e^{-2\lambda s}R^2$. The
$V_2\to H$ error bound yields
\[
 \norm{I_{1,\mathrm{late}}}_H\le ChR^2\int_{t/2}^t(t-s)_\wedge^{-1/2}
 s_\wedge^{-3}e^{-\nu_1(t-s)-2\lambda s}\,ds
 \le Cht_\wedge^{-5/2}e^{-\nu t}R^2 .
\]

\emph{Early part of $I_2$.} Put $\alpha_h(s)=\min\{1,hs_\wedge^{-3/2}\}$
and $\beta(s)=e^{\nu_1s}(|\Phi(s)|_1+|\Phi_h(s)|_1)$. By the contraction property and the $H\to H$ error bound, $\norm{\Delta_h(s)}_H\le CRe^{-\nu_1s}\alpha_h(s)$. Moreover, by Agmon's inequality, $\norm{\Phi+\Phi_h}_\infty\le CR^{1/2}e^{-\nu_1s}\beta(s)^{1/2}$. Since
$\ip{F(f,g)}w=\frac12\int_0^Lfg\,\overline{(\rho w)'}\dx$,
$\norm{F(\Delta_h,\Phi+\Phi_h)}_{V_1^*}\le C\norm{\Delta_h}_H
\norm{\Phi+\Phi_h}_\infty$. By \eqref{eq:weaksmoothing}, $t-s\ge t/2$ and
H\"older's inequality with exponents $4/3$ and $4$,
\begin{align*}
 \norm{I_{2,\mathrm{early}}}_H
 &\le Ct_\wedge^{-3/4}e^{-\nu t}R^{3/2}
 \Big(\int_0^\infty e^{-\frac43(2\nu_1-\nu)s}\alpha_h(s)^{4/3}ds\Big)^{3/4}
 \Big(\int_0^\infty\beta(s)^2ds\Big)^{1/4}\\
 &\le Ch^{1/2}t_\wedge^{-3/4}e^{-\nu t}R^2 .
\end{align*}
Indeed, for fixed $c_0>0$ and $h\le1$,
\[
 \int_0^\infty e^{-c_0s}\alpha_h(s)^{4/3}\,ds
 \le h^{2/3}+h^{4/3}\int_{h^{2/3}}^1s^{-2}\,ds+C_{c_0}h^{4/3}
 \le C_{c_0}h^{2/3}.
\]
This is the only step that loses half a power of $h$.

\emph{Late part of $I_2$.} By \eqref{eq:cpH1} and Agmon's inequality,
$\norm{\Phi+\Phi_h}_\infty\le Cs_\wedge^{-3/8}e^{-\lambda s/2}R$. With
the $H\to H$ error bound and \eqref{eq:weaksmoothing},
\[
 \norm{I_{2,\mathrm{late}}}_H\le ChR^2\int_{t/2}^t(t-s)_\wedge^{-3/4}
 s_\wedge^{-15/8}e^{-\lambda(t-s)-(\nu_1+\lambda/2)s}\,ds
 \le Cht_\wedge^{-13/8}e^{-\nu t}R^2 .
\]
Adding the four bounds and using $h\le h^{1/2}$, we obtain
\eqref{eq:energyquadratic}.
\end{proof}

The rate $h^{1/2}$ holds for the quadratic coefficient on the whole
energy ball. Proposition~\ref{thm:energyfull} below bounds the complete
ordering error and its cubic remainder, at the lower rate $h^{1/6}$. Write
\[
 \norm w_{\Z_\nu}=\Big(\int_0^\infty e^{4\nu t}\norm{w(t)}_H^4\,dt\Big)^{1/4}.
\]
All trajectory spaces are built on $H$.
Discrete trajectories are included via $X_h\subset H_0^1\subset H$.

\begin{lemma}[Mixed time norms]\label{lem:energymixed}
For complex trajectories in $\X_0$,
\begin{align}
 \norm w_{\Z_\nu}+\norm w_{L^4_tL^\infty_x}
   &\le C\E_\nu(w),\label{eq:mixedembed}\\
 \norm{e^{\nu t}\N(w)}_{L^{4/3}_tH}
   &\le C\E_\nu(w)^2,\label{eq:mixedsource}\\
 \norm{e^{\nu t}F(f,w)}_{L^2_tV_1^*}
   &\le C\norm f_{\Z_\nu}\E_\nu(w).\label{eq:mixedfeedback}
\end{align}
In the last estimate $f\in\Z_\nu$ suffices. Moreover,
\begin{align}
 \norm{L_h(\cdot)f}_{\Z_\nu}&\le Ch^{1/6}\norm f_H,\label{eq:Jinitial}\\
 \norm{L_h*g}_{\Z_\nu}&\le Ch^{1/3}\norm{e^{\nu t}g}_{L^{4/3}_tH},
 \label{eq:JconvZ}\\
 \norm{L_h*g}_{\Y_\nu(H)}&\le Ch^{1/6}\norm{e^{\nu t}g}_{L^{4/3}_tH}.
 \label{eq:JconvY}
\end{align}
\end{lemma}
\begin{proof}
Let $M_w=\sup_te^{\nu t}\norm{w(t)}_H$ and
$J_w^2=\int_0^\infty e^{2\nu t}|w(t)|_1^2\,dt$. By Poincar\'e's inequality, we have
\[
 \int_0^\infty e^{4\nu t}\norm{w}_H^4\,dt\le
 M_w^2\int_0^\infty e^{2\nu t}\norm w_H^2\,dt\le p^2M_w^2J_w^2 .
\]
By Agmon's inequality, $\norm{w(t)}_\infty^4\le C\norm{w(t)}_H^2|w(t)|_1^2$.
Its unweighted time integral is at most $CM_w^2J_w^2$. For almost every
$t$, $\norm{\N(w(t))}_H\le C\norm{w(t)}_H^{1/2}|w(t)|_1^{3/2}$, hence
\[
 \int_0^\infty e^{4\nu t/3}\norm{\N(w(t))}_H^{4/3}\,dt
 \le CM_w^{2/3}\int_0^\infty e^{2\nu t/3}|w(t)|_1^2\,dt\le CM_w^{2/3}J_w^2 .
\]
Taking the $3/4$ power, we obtain \eqref{eq:mixedsource}. By the weighted weak pairing, $\norm{F(f,w)}_{V_1^*}\le C\norm f_H\norm w_\infty$.
By H\"older's inequality in time, with $\frac12=\frac14+\frac14$, we obtain
\eqref{eq:mixedfeedback}. The whole exponential weight is placed on $f$.
No inverse estimate is used.

For the last three estimates, we note that $e^{\nu t}\norm{L_h(t)f}_H\le
\Lambda_h(t)\norm f_H$. This yields \eqref{eq:Jinitial} by
\eqref{eq:kernels} with $p_0=4$. For the convolutions, we extend the kernel
$\Lambda_h$ and the function $e^{\nu t}\norm{g(t)}_H$ by zero to negative
times. By Young's inequality with $1+\frac14=\frac12+\frac34$, we obtain \eqref{eq:JconvZ}, and with $1=\frac14+\frac34$ we obtain \eqref{eq:JconvY}.
\end{proof}

\begin{proposition}[The complete ordering error on an energy ball]
\label{thm:energyfull}
There are $r,C>0$, independent of $h$, such that every real or complex
$u_0\in H$ with $\varepsilon=\norm{u_0}_H\le r$ satisfies
\begin{align}
 \norm{u-u_h}_{\Z_\nu}&\le Ch^{1/6}\varepsilon,
       \label{eq:energyfullZ}\\
 \norm{\Cd_h(u_0)}_{\Y_\nu(H)}&\le Ch^{1/6}\varepsilon^2,
       \label{eq:energyfulldefect}\\
 \norm{\Cd_h(u_0)-D_h(u_0)}_{\Y_\nu(H)}&\le Ch^{1/6}\varepsilon^3.
       \label{eq:energyfulltail}
\end{align}
The last bound controls the sum of all Taylor terms of degree at least
three. No higher norm of $u_0$ and no $h$-dependent smallness condition
is required.
\end{proposition}
\begin{proof}
By Theorem~\ref{thm:continuous} and
Proposition~\ref{prop:uniformanalytic} there exist a common complex ball
$B_H(0,r_{\rm an})$ and a constant $M_0$, independent of $h$, with
$\E_\nu(u)+\E_\nu(u_h)\le M_0\norm{u_0}_H$. For the discrete map we first
apply the complex-linear contraction $\Pi_h$ to $u_0$ and then embed the
trajectory in the common space $\X_0$. Both maps are holomorphic with
values in this same Banach space. By \eqref{eq:disforced}, we have
$\E_\nu(\QQ_h(y_1,y_2))\le C_F\norm{e^{\nu t}F(y_1,y_2)}_{L^2_tV_1^*}$
whenever the right side is finite. Here $F_h=\Pi_hF$ acts on a weak
source by restriction to $X_h$. Its weak norm does not increase.

It follows from the identities \eqref{eq:splitting} and \eqref{eq:Cdmild} that
\begin{equation}\label{eq:fullid}
 e=-L_h(\cdot)u_0-L_h*\N(u)+2\QQ_h(e,m).
\end{equation}
These are identities between energy solutions. In particular,
$e,m\in\X_0$ and $e\in\Z_\nu$ before any absorption. The source estimate \eqref{eq:mixedsource} yields local $L^1_tH$ integrability. Hence the
mild integrals and the weak forced solutions coincide. By
Lemma~\ref{lem:energymixed},
\[
 \norm e_{\Z_\nu}\le C_1h^{1/6}\varepsilon+C_2h^{1/3}\varepsilon^2
 +C_3\varepsilon\norm e_{\Z_\nu}.
\]
Choose $0<R_0<r_{\rm an}$ with $C_3R_0\le1/2$. For complex
$u_0$ with $\norm{u_0}_H\le R_0$, absorption and $h^{1/3}\le h^{1/6}$ yield
\eqref{eq:energyfullZ}. Then, by \eqref{eq:JconvY} and the forced estimate, applied to \eqref{eq:Cdmild}, we obtain
\[
 \norm{\Cd_h(u_0)}_{\Y_\nu(H)}
 \le Ch^{1/6}\E_\nu(u)^2+C\norm e_{\Z_\nu}\E_\nu(m)
 \le \bar Ch^{1/6}\varepsilon^2 .
\]
All trajectories are $H$-continuous. Hence these estimates hold at every
time. At zero the ordering error vanishes.

It remains to estimate the Taylor tail. For a complex unit vector
$\phi\in H$, the function $z\mapsto\Cd_h(z\phi)$ is holomorphic with
values in $\Y_\nu(H)$ on a neighborhood of $|z|\le R_0$. Its constant
and linear coefficients vanish, and its quadratic coefficient is
$D_h(\phi)$. Write $\Cd_h(z\phi)=\sum_{n\ge2}z^na_{n,h}(\phi)$. The
preceding estimate holds on the circle $|z|=R_0$. Hence, by the Banach-valued Cauchy formula,
$\norm{a_{n,h}(\phi)}_{\Y_\nu(H)}\le\bar Ch^{1/6}R_0^{2-n}$ for $n\ge2$. Put
$r=R_0/2$. For $u_0=\varepsilon\phi$ with $0<\varepsilon\le r$,
\[
 \norm{\Cd_h(u_0)-D_h(u_0)}_{\Y_\nu(H)}
 \le\bar Ch^{1/6}R_0^2\sum_{n\ge3}(\varepsilon/R_0)^n
 \le\frac{2\bar C}{R_0}h^{1/6}\varepsilon^3 .
\]
The constants and radii are independent of the direction, of the mesh and
of the number of Taylor terms. The proof is complete.
\end{proof}

The factor $h^{1/6}$ is the $L^4$ size of an initial layer of width $h^{2/3}$.
The proof uses first-order consistency only for strong sources, in
$L^{4/3}_tH$, and only stability for weak sources.

\emph{Proof of Theorems~\ref{thm:linear} and~\ref{thm:ordering}.} The
four solutions exist globally by Theorem~\ref{thm:continuous} and
Proposition~\ref{prop:nonlinstab}. The identity \eqref{eq:splitting}
holds by construction. The estimate \eqref{eq:sn-linear} follows from
Proposition~\ref{thm:positive} and the contraction property. The
estimate \eqref{eq:sn-full} follows from \eqref{eq:splitting},
\eqref{eq:sn-linear} and \eqref{eq:sn-ordering}. Since the bounds in
\eqref{eq:sn-linear} and \eqref{eq:sn-full} are nonincreasing in $t$,
uniform convergence on $[t_0,\infty)$ follows. For $u-v$, the bound
\eqref{eq:sn-horizontal} follows from the continuous Taylor remainder
\eqref{eq:taylor} with $m=1$. For $u_h-v_h$,
Proposition~\ref{prop:uniformanalytic} yields this bound. The estimate
\eqref{eq:sn-ordering}, including $t=0$, is \eqref{eq:energyfulldefect}.
The Taylor tail is \eqref{eq:energyfulltail}, and the bound for $D_h$ is
Proposition~\ref{prop:energyquadratic}. The radius was reduced finitely
many times, each time to a value depending only on $L$, $\kappa$ and the
fixed rates.\hfill$\square$

\subsection{\texorpdfstring{Data in $D(A)$}{Data in D(A)}}
\label{sec:convproof}
In this subsection, we prove Theorem~\ref{thm:DA}. Part~(i) is the
second estimate in \eqref{eq:H4integrated}, since
$v-v_h=-L_h(\cdot)u_0$. The proof of parts~(ii) and~(iii) is divided into
four steps: the convergence of the linear maps, the convergence in the
energy norm, the rate for data in $D(A)$, and the ordering error for
data in $D(A)$.

Write $\mathcal Kg(t)=\int_0^tT(t-s)g(s)\,ds$ and
$\mathcal K_hg(t)=\int_0^tT_h(t-s)\Pi_hg(s)\,ds$ for the forced solutions
with zero initial value. For weak sources they are the energy solutions
of Sections~\ref{sec:analyticproof} and~\ref{sec:stability}. Let
$L^2_\nu(V_1^*)$ be the space of sources with
$\norm g_{L^2_\nu(V_1^*)}=\norm{e^{\nu t}g}_{L^2_tV_1^*}<\infty$. By
\eqref{eq:forced} and \eqref{eq:disforced}, and by the weak product bound
of Section~\ref{sec:analyticproof},
\begin{equation}\label{eq:Kbounds}
 \E_\nu(\mathcal Kg)+\E_\nu(\mathcal K_hg)\le C\norm g_{L^2_\nu(V_1^*)},
 \qquad
 \norm{F(w,\tilde w)}_{L^2_\nu(V_1^*)}\le C_B\E_\nu(w)\E_\nu(\tilde w).
\end{equation}

\emph{Step 1: linear maps.} Both linear flows are bounded from $H$ into
$\X_0$, and $\E_\nu(L_h(\cdot)g)\le Ch^2\norm{Ag}_H$ for $g\in D(A)$, by
Lemma~\ref{lem:DAenergy}. Since $D(A)$ is dense in $H$,
$\E_\nu(v-v_h)=\E_\nu(L_h(\cdot)u_0)\to0$ for every $u_0\in H$. The same
argument applies to forced solutions. For $G\in C_c^\infty((0,\infty);D(A))$,
$(\mathcal K-\mathcal K_h)G(t)=-\int_0^tL_h(t-s)G(s)\,ds$. A time shift
by $s$ multiplies the weighted norm by $e^{\nu s}$. Hence, by Minkowski's inequality and Lemma~\ref{lem:DAenergy}, we obtain
\[
 \E_\nu((\mathcal K-\mathcal K_h)G)
 \le Ch^2\int_0^\infty e^{\nu s}\norm{AG(s)}_H\,ds .
\]
Such sources are dense in $L^2_\nu(V_1^*)$. Indeed, we truncate and mollify in time
and then approximate in space, since $D(A)$ is dense in $H$ and $H$ is
dense in $V_1^*$ (by reflexivity of $V_1$, since $V_1$ is dense in $H$).
By the uniform bound \eqref{eq:Kbounds},
\begin{equation}\label{eq:forcedstrong}
 \E_\nu((\mathcal K-\mathcal K_h)G)\longrightarrow0\qquad
 \text{for each fixed }G\in L^2_\nu(V_1^*).
\end{equation}

\emph{Step 2: convergence in the energy norm.} By
Theorem~\ref{thm:continuous} and Proposition~\ref{prop:uniformanalytic},
$\E_\nu(u)+\E_\nu(u_h)\le M_0\norm{u_0}_H$ for $\norm{u_0}_H\le r$.
Subtracting the two mild equations and using $\N(u)-\N(u_h)=2F(e,m)$, we obtain
\[
 e=-L_h(\cdot)u_0+(\mathcal K-\mathcal K_h)\N(u)+2\mathcal K_hF(e,m).
\]
By \eqref{eq:Kbounds}, the last term has energy norm at most
$CC_BM_0r\,\E_\nu(e)\le\frac12\E_\nu(e)$ for a small fixed $r$. Since
$\E_\nu(e)$ is finite,
\[
 \E_\nu(e)\le2\E_\nu(L_h(\cdot)u_0)
 +2\E_\nu((\mathcal K-\mathcal K_h)\N(u))\longrightarrow0,
\]
by Step 1 and \eqref{eq:forcedstrong} with the fixed source
$\N(u)\in L^2_\nu(V_1^*)$. No error
estimate for the nonlinear problem has been used.

\emph{Step 3: Theorem~\ref{thm:DA}(ii).} Let $g\in D(A)$ be real with
$\norm g_H\le r$, and let $u_0=g$. Its $H^3$ norm may be large. The
constants below do not depend on it. We first bound the discrete
solutions, which solve ordinary differential equations, and then let
$h\to0$.

The right side of \eqref{eq:Wscheme} is a quadratic polynomial on $X_h$.
Hence $u_h$ is smooth in time, and $q_h=\dot u_h$ solves the linear
equation
\[
 \dot q_h=A_hq_h+2F_h(u_h,q_h),\qquad
 q_h(0)=A_h\Pi_hg+\Pi_h\N(\Pi_hg).
\]
By \eqref{eq:graphproj}, $\norm{A_h\Pi_hg}_H\le C\norm g_{H^3}$. By
\eqref{eq:prodenergy} and \eqref{eq:disenergy},
$\norm{\N(\Pi_hg)}_H\le C|\Pi_hg|_1^2\le C\norm g_H\norm{A_h\Pi_hg}_H$.
Hence $\norm{q_h(0)}_H\le C\norm g_{H^3}$. On $\X_0$ consider the affine
map $\Phi(q)=T_h(\cdot)q_h(0)+2\mathcal K_hF(u_h,q)$. By \eqref{eq:Kbounds}
and Step~2, its linear part has norm at most $2CC_BM_0\norm g_H$, which
is at most $1/2$ for a small fixed $r$. Thus $\Phi$ is a contraction.
Its fixed point solves the same linear equation with the same initial
value as $q_h$. By uniqueness for linear ordinary differential equations,
the fixed point is $q_h$. Hence $\E_\nu(\dot u_h)\le C\norm g_{H^3}$.
The scheme also gives $B_hu_h=\Pi_h\N(u_h)-q_h$. As above,
$\norm{\N(u_h)}_H\le C\norm{u_h}_H\norm{B_hu_h}_H$, and
$\norm{u_h(t)}_H\le\norm g_H$ by Proposition~\ref{prop:nonlinstab}. For
small $r$ this term is absorbed, and
$\norm{A_hu_h(t)}_H\le2\norm{q_h(t)}_H\le Ce^{-\nu t}\norm g_{H^3}$.

Now let $h\to0$. By Step~2, $\E_\nu(u-u_h)\to0$. By \eqref{eq:graph},
$\norm{u_h(t)}_{H^2}\le Ce^{-\nu t}\norm g_{H^3}$, so $u_h(t)\to u(t)$
weakly in $H^2$. For $\varphi\in C_c^\infty(0,L)$ we have
$\Pi_h\varphi\to\varphi$ in $H^1$ by \eqref{eq:approx}, and
\[
 |b(u_h(t),\Pi_h\varphi)|=|\ip{A_hu_h(t)}{\Pi_h\varphi}|
 \le Ce^{-\nu t}\norm g_{H^3}\norm\varphi_H .
\]
In the limit, $|\ip{u(t)'''}\varphi|\le Ce^{-\nu t}\norm g_{H^3}\norm\varphi_H$.
Hence $u(t)\in D(A)$ and $\norm{u(t)}_{H^3}\le Ce^{-\nu t}\norm g_{H^3}$.
After extraction of a subsequence, $e^{\nu t}\dot u_h$ converges weak-star
in $L^\infty(0,\infty;H)$ and weakly in $L^2(0,\infty;H_0^1)$ to
$e^{\nu t}q$. Passing to the limit in $u_h(t)-u_h(s)=\int_s^t\dot u_h\,dr$
gives $q=\dot u$. By weak lower semicontinuity,
$\E_\nu(\dot u)\le C\norm g_{H^3}$, where the supremum is at first an
essential supremum. Since $Au$ is bounded in $H$, the mild equation gives
$\dot u=Au+\N(u)$ almost everywhere. Moreover
$\dot u\in L^2_{\rm loc}(H_0^1)$, so $\N(u)$ is absolutely continuous in
$H$ with derivative $2F(u,\dot u)$, by \eqref{eq:prodenergy}.
Differentiating the mild equation, we see that $\dot u$ is the energy
solution of $\dot q=Aq+2F(u,q)$ with $q(0)=Ag+\N(g)$. In particular
$\dot u$ is continuous in $H$, and the essential supremum is a supremum.
Thus $u$ is a global strong solution, and $u'''=Bu=\N(u)-\dot u$. Since
$\norm{\N(u)'}_2\le C|u|_1\norm u_{H^3}$,
$\norm{u(t)}_{H^3}\le Ce^{-\nu t}\norm g_{H^3}$ and
$\int_0^\infty|u|_1^2dt\le C\norm g_H^2$, we have
\[
 \norm{e^{\nu t}\N(u)'}_{L^2_tL^2_x}\le C\norm g_{H^3}\norm g_H .
\]
With $\norm w_{H^4}\le C(\norm w_H+\norm{w''''}_2)$, differentiating the
identity once in $x$ yields the integrated $H^4$ bound, linear in
$\norm g_{H^3}$, as in the proof of Lemma~\ref{lem:DAenergy}. Hence
\begin{equation}\label{eq:nonlinH4}
 \sup_{t\ge0}e^{\nu t}\norm{u(t)}_{H^3}
 +\norm{e^{\nu t}u}_{L^2_tH^4}\le C\norm g_{H^3}.
\end{equation}
Now write $\zeta=u-\Pi_hu$ and $y_h=\Pi_hu-u_h$. Then
$\dot y_h=A_hy_h+2F_h(m,y_h)+d_h(u)+2F_h(m,\zeta)$, $y_h(0)=0$. By
\eqref{eq:approx} and \eqref{eq:nonlinH4},
$\E_\nu(\zeta)\le Ch^3\norm g_{H^3}$. The residual bound \eqref{eq:residual} with $s=4$ gives
$\norm{e^{\nu t}d_h(u)}_{L^2_t(-1,h)}\le Ch^2\norm g_{H^3}$, and
\eqref{eq:prodenergy} gives
$\norm{e^{\nu t}F_h(m,\zeta)}_{L^2_t(-1,h)}\le
C\sup_te^{\nu t}\norm\zeta_H\norm{|m|_1}_{L^2_t}\le Ch^3\norm g_{H^3}$.
It then follows from \eqref{eq:disforced}, \eqref{eq:Kbounds} and Step~2 that
$\E_\nu(y_h)\le Ch^2\norm g_{H^3}+2CC_BM_0\norm g_H\E_\nu(y_h)$. Here
$\E_\nu(y_h)$ is finite, since $\Pi_h$ is bounded on $V_1$. For small
$r$ the last term is absorbed. With the bound for $\zeta$, we obtain
\begin{equation}\label{eq:DAnonlinear}
 \E_\nu(u-u_h)\le Ch^2\norm g_{H^3}.
\end{equation}
Together with Lemma~\ref{lem:DAenergy} for $v-v_h$, this proves
Theorem~\ref{thm:DA}(ii).

\emph{Step 4: Theorem~\ref{thm:DA}(iii).} Let $u_0\in D(A)$ be real with
$\varepsilon=\norm{u_0}_H\le r$, where $r$ is the radius of Step~3. By
\eqref{eq:Cdmild},
\[
 \Cd_h=-L_h*\N(u)+2\QQ_h(e,m).
\]
We use three facts. By Step~3, $u(t)\in D(A)$ for $t\ge0$, and
\eqref{eq:nonlinH4} holds with $g=u_0$. Next,
$\E_\nu(u)+\E_\nu(u_h)\le M_0\varepsilon$, by Step~2. Finally,
$\E_\nu(e)\le Ch^2\norm{u_0}_{H^3}$, by \eqref{eq:DAnonlinear}.

First, by the interpolation inequality
$\norm w_{H^2}^2\le C\norm w_{H^1}\norm w_{H^3}$, by Poincar\'e's
inequality and by the Cauchy--Schwarz inequality in time,
\[
 \int_0^\infty\norm{u}_{H^2}^2\,dt
 \le C\sup_{t\ge0}e^{\nu t}\norm{u}_{H^3}\int_0^\infty e^{-\nu t}\norm{u'}_2\,dt
 \le C\norm{u_0}_{H^3}\,\E_\nu(u)\le C\norm{u_0}_{H^3}\,\varepsilon .
\]

Second, the source stays in $D(A)$. Let $w=u(s)\in H^4\cap D(A)$; by
\eqref{eq:nonlinH4} this holds for almost every $s$. Then
$G=\N(w)=-ww'\in H^3$, and $G$ vanishes at $0$ and $L$. Moreover
$G'=-(w')^2-ww''$, so
\[
 G'(L)=-w'(L)^2-w(L)w''(L)=0,
\]
because $w(L)=w'(L)=0$. Hence $G\in D(A)$. From
$G'''=-(ww''''+4w'w'''+3(w'')^2)$ and $\norm{w^{(k)}}_\infty\le C\norm w_{H^{k+1}}$
we obtain $\norm{AG}_H\le C\norm w_{H^2}\norm w_{H^4}$. By the
Cauchy--Schwarz inequality, \eqref{eq:nonlinH4} and the first step,
\begin{equation}\label{eq:sourceint}
 \int_0^\infty e^{\nu s}\norm{A\N(u(s))}_H\,ds
 \le C\norm{u_0}_{H^3}^{3/2}\varepsilon^{1/2}.
\end{equation}
The map $s\mapsto\N(u(s))$ is strongly measurable into $D(A)$,
since $u\in L^2_{\rm loc}(0,\infty;H^4)$.

Third, $\E_\nu(L_h(\cdot)g)\le Ch^2\norm{Ag}_H$ for $g\in D(A)$, by
\eqref{eq:H4integrated}. The function $t\mapsto L_h(t-s)g\,\mathbf 1_{t\ge s}$
is the shift of $L_h(\cdot)g$ by $s$, and the shift multiplies both parts
of $\E_\nu$ by $e^{\nu s}$. By Minkowski's inequality, applied to
each part, and by \eqref{eq:sourceint},
\[
 \E_\nu(L_h*\N(u))
 \le\int_0^\infty e^{\nu s}\E_\nu\big(L_h(\cdot)\N(u(s))\big)\,ds
 \le Ch^2\norm{u_0}_{H^3}^{3/2}\varepsilon^{1/2}.
\]
This is the argument of Step~1.

Fourth, by \eqref{eq:Kbounds} and the three facts,
$\E_\nu(\QQ_h(e,m))\le C\E_\nu(e)\E_\nu(m)
\le Ch^2\norm{u_0}_{H^3}\varepsilon$. Since $\varepsilon\le C\norm{u_0}_{H^3}$,
this is at most $Ch^2\norm{u_0}_{H^3}^{3/2}\varepsilon^{1/2}$. Adding the
two terms and using $\norm{u_0}_{H^3}\le C\norm{Au_0}_H$, we obtain
the first estimate in Theorem~\ref{thm:DA}(iii). The bound at each time
follows from the supremum in $\E_\nu$. For the last assertion, the
definition \eqref{eq:defect} of the ordering error on the mesh of size
$\hat h$ gives $u-\tilde u=(v-v_h)+\Cd_{\hat h}$. By (i),
$e^{\nu t}\norm{(v-v_h)(t)}_H\le\E_\nu(v-v_h)\le Ch^2\norm{Au_0}_H$. The
second term is bounded by the estimate at each time on the mesh of size
$\hat h$. This completes the proof of
Theorem~\ref{thm:DA}.\hfill$\square$

\section{Optimal design of the two-mesh computation}\label{sec:design}
\setcounter{equation}{0}

In this section, we design the two-mesh computation at minimal cost. We
use the W-scheme and uniform meshes of sizes $h=L/N$ and $\hat h=L/\hat N$
with $\hat N\le N$ and $\hat h\le\min\{1,L\}$, as in
Section~\ref{sec:scheme}. They need not be nested. A nonlinear simulation
costs more than a linear one, since every time step requires an
iteration. It is therefore natural to run the two on different meshes.
Recall the two-mesh sum
$\tilde u=v_h+(u_{\hat h}-v_{\hat h})$ and its error
$u-\tilde u=(v-v_h)+\Cd_{\hat h}$ from Theorem~\ref{thm:DA}(iii).

\subsection{The accuracy constraint}\label{sec:design-constraint}
Fix $0<\nu<\lambda$, and fix $0<t_0\le1$, $\eta>0$
and $0<\varepsilon\le r$, where $r$ is the smallest of the radii
in Theorems~\ref{thm:linear}, \ref{thm:ordering} and~\ref{thm:DA}. We require
\begin{equation}\label{eq:accuracy}
 \norm{(u-\tilde u)(t)}_H\le\eta\norm{u_0}_H
 \quad\text{for all }t\ge t_0\text{ and all real }u_0\text{ with }\norm{u_0}_H\le\varepsilon .
\end{equation}
Let $t\ge t_0$ and $\norm{u_0}_H\le\varepsilon$. Then $e^{-\nu t}\le1$ and
$t_\wedge=\min\{1,t\}\ge t_0$. By Theorem~\ref{thm:linear}, we obtain
$\norm{(v-v_h)(t)}_H\le C_1t_0^{-3/2}h\norm{u_0}_H$.
By Theorem~\ref{thm:ordering}, the ordering error is the discretization
error $D_{\hat h}(u_0)$ of the quadratic Taylor coefficient up to a cubic
remainder, with
\[
 \norm{D_{\hat h}(u_0)(t)}_H\le C_2t_0^{-5/2}\hat h^{1/2}\norm{u_0}_H^2,\qquad
 \norm{(\Cd_{\hat h}-D_{\hat h}(u_0))(t)}_H\le C_3\hat h^{1/6}\norm{u_0}_H^3 .
\]
The constants $C_1,C_2,C_3$ depend only on $L$ and $\nu$. They do not
depend on $h$, $\hat h$, $t$, $t_0$, $\varepsilon$, $\eta$ and $u_0$.
Hence \eqref{eq:accuracy} holds if
\begin{equation}\label{eq:PhiPsi}
 \Phi(h)+\Psi(\hat h)\le\eta,\qquad
 \Phi(h)=C_1t_0^{-3/2}h,\qquad
 \Psi(\hat h)=C_2t_0^{-5/2}\varepsilon\,\hat h^{1/2}+C_3\varepsilon^2\hat h^{1/6}.
\end{equation}
The same bounds with $\hat h=h$ apply to the full error
$u-u_h=(v-v_h)+\Cd_h$ of the nonlinear solution on one mesh. So the bound
$\Phi(h)+\Psi(h)\le\eta$ implies \eqref{eq:accuracy} for $u_h$. The factor
$\varepsilon$ in $\Psi$ comes from the mesh-uniform product bound for the
ordering error. It separates the two meshes: the fine mesh serves the
linear part, and the coarse mesh serves a term that is small with the
datum.

\subsection{The cost function}\label{sec:design-cost}
We measure cost by computing time on one processor thread. On a mesh
with $N$ elements the scheme has $n=2N-1=2L/h-1$ unknowns. We assume
that one time step of a linear run takes the time $c_{\rm L}+d_{\rm L}n$,
and one time step of a nonlinear run the time $c_{\rm N}+d_{\rm N}n$,
with constants $c_{\rm L},d_{\rm L},c_{\rm N},d_{\rm N}>0$. The part
proportional to $n$ is the work on all unknowns, such as the solution of
the banded linear systems. The fixed part does not depend on the mesh. A
nonlinear step solves a banded system in each of its iterations, so it
takes longer than a linear step. We assume $d_{\rm N}>d_{\rm L}$. We also
assume that every run makes the same number $N_t$ of time steps. Thus, on
one processor thread, the times of a linear run and of a nonlinear run on
a mesh of size $h$ are
\[
 J_{\rm L}(h)=N_t(c_{\rm L}+d_{\rm L}n),\qquad J_{\rm N}(h)=N_t(c_{\rm N}+d_{\rm N}n).
\]
The cost function is the serial cost of $\tilde u$, the computing time of
its three runs one after another:
\[
 J(h,\hat h)=J_{\rm L}(h)+J_{\rm L}(\hat h)+J_{\rm N}(\hat h)
 =J_0+\frac{K_1}h+\frac{K_2}{\hat h},\qquad
 K_1=2LN_td_{\rm L},\quad K_2=2LN_t(d_{\rm L}+d_{\rm N}),
\]
where $J_0$ does not depend on the meshes. The cost of $u_h$ is
$J_{\rm N}(h)=J_0'+K/h$ with $K=2LN_td_{\rm N}$. We treat the mesh sizes as
continuous variables and consider the nonlinear programming model
\begin{equation}\label{eq:design}
 \text{minimize }\ \frac{K_1}h+\frac{K_2}{\hat h}
 \quad\text{subject to}\quad \Phi(h)+\Psi(\hat h)\le\eta,\quad 0<h\le\hat h\le\bar h ,
\end{equation}
where $\bar h\le\min\{1,L\}$ is the largest allowed mesh size. We assume
\begin{equation}\label{eq:hbar}
 \Phi(\bar h)+\Psi(\bar h)>\eta .
\end{equation}
Otherwise $h=\hat h=\bar h$ satisfies the accuracy bound
and minimizes the cost. A pair
$(h,\hat h)$ that satisfies the constraints is called a feasible solution.

\subsection{The optimal design}\label{sec:design-opt}
Write the constraints of \eqref{eq:design} as
\[
 g_0=\Phi(h)+\Psi(\hat h)-\eta\le0,\qquad g_1=h-\hat h\le0,\qquad
 g_2=\hat h-\bar h\le0 .
\]
The cost tends to infinity as $h\to0^+$. The constraint $h>0$ is open and
needs no multiplier. The Lagrange function is
\[
 \mathcal L(h,\hat h,\xi_0,\xi_1,\xi_2)=\frac{K_1}h+\frac{K_2}{\hat h}
 +\xi_0g_0+\xi_1g_1+\xi_2g_2 .
\]
A feasible solution $(h,\hat h)$ satisfies the Karush--Kuhn--Tucker
conditions~\cite[Section~4.2]{BSS} if there exist $\xi_0,\xi_1,\xi_2\ge0$ such
that $\nabla_{(h,\hat h)}\mathcal L=0$ and $\xi_ig_i=0$ for $i=0,1,2$, that
is,
\begin{equation}\label{eq:KKT}
 \begin{gathered}
 -\frac{K_1}{h^2}+\xi_0\Phi'(h)+\xi_1=0,\qquad
 -\frac{K_2}{\hat h^2}+\xi_0\Psi'(\hat h)-\xi_1+\xi_2=0,\\
 \xi_0\bigl(\Phi(h)+\Psi(\hat h)-\eta\bigr)=0,\qquad
 \xi_1(h-\hat h)=0,\qquad \xi_2(\hat h-\bar h)=0 .
 \end{gathered}
\end{equation}
The feasible set of \eqref{eq:design} is nonconvex in $(h,\hat h)$.
By \eqref{eq:hbar}, the accuracy boundary meets the diagonal below $\bar h$
and extends into $0<h<\hat h<\bar h$. On this arc, $h$ is a strictly convex
function of $\hat h$, since $\Phi$ is linear and increasing and $\Psi$ is
strictly concave. The midpoint of two distinct points on the arc lies
above the boundary and violates the accuracy constraint.
The problem is convex in the logarithms of the mesh sizes.

\begin{theorem}[Optimal design]\label{prop:design}
Let $\Phi$ and $\Psi$ be defined as in \eqref{eq:PhiPsi}, and let
\eqref{eq:hbar} hold. Then problem \eqref{eq:design} has a
minimizer $(h^*,\hat h^*)$. It satisfies the Karush--Kuhn--Tucker system
\eqref{eq:KKT}, and this system has exactly one solution. Moreover, $\Phi(h^*)+\Psi(\hat h^*)=\eta$.
\end{theorem}

Let $0<h_1\le\bar h$ be the solution of $\Phi(h_1)+\Psi(h_1)=\eta$. Since
$\Phi+\Psi$ is increasing, $h_1$ is the coarsest single mesh for which the
bound $\Phi(h)+\Psi(h)\le\eta$ holds.

\begin{corollary}[Structure of the optimum]\label{cor:design-structure}
Under the assumptions of Theorem~\ref{prop:design}, the following hold.
\begin{enumerate}
\item[(a)] $h^*\le h_1\le\hat h^*$.
\item[(b)] $h^*=\hat h^*=h_1$ if and only if
$K_2\Phi'(h_1)\le K_1\Psi'(h_1)$.
\item[(c)] If $h^*<\hat h^*<\bar h$, then
\[
 \frac{K_1}{(h^*)^2\Phi'(h^*)}=\frac{K_2}{(\hat h^*)^2\Psi'(\hat h^*)} .
\]
Moreover, the optimal mesh ratio is
\[
 \frac{\hat h^*}{h^*}=\frac{C_1t_0^{-3/2}x^6}{\eta-\Psi(x^6)},
\]
where $x=(\hat h^*)^{1/6}$ is the only zero of the polynomial
\[
 P(x)=\frac{K_1C_1}{K_2t_0^{3/2}}\Bigl(\frac{C_2\varepsilon}{2t_0^{5/2}}x^9+\frac{C_3\varepsilon^2}{6}x^7\Bigr)
 -\Bigl(\eta-\frac{C_2\varepsilon}{t_0^{5/2}}x^3-C_3\varepsilon^2x\Bigr)^2
\]
with $x>0$ and $\Psi(x^6)<\eta$.
\item[(d)] If $h^*<\hat h^*<\bar h$ for all $\eta$ in an interval, then
on this interval $h^*$ and $\hat h^*$ increase strictly with respect to $\eta$, and so
does the mesh ratio
\[
 \frac{\hat h^*}{h^*}=\Bigl(\frac{C_1t_0^{-3/2}}{\Psi'(\hat h^*)}\cdot\frac{K_2}{K_1}\Bigr)^{1/2}.
\]
\end{enumerate}
\end{corollary}

\begin{proof}[Proof of Theorem~\ref{prop:design}]
For small $h=\hat h$ the constraint holds, so the feasible set is not
empty. The feasible sublevel set determined by any feasible pair is
compact. On this set $h$ is bounded away from zero by the cost bound,
and $h\le\hat h\le\bar h$. Hence a minimizer exists. Let
$(h,\hat h)$ be a minimizer, and suppose $\Phi(h)+\Psi(\hat h)<\eta$. If
$\hat h<\bar h$, a slightly larger $\hat h$ is feasible and cheaper. So
$\hat h=\bar h$. In the same way $h=\hat h$. Then
$\Phi(\bar h)+\Psi(\bar h)<\eta$, a contradiction. So every minimizer
satisfies $\Phi(h)+\Psi(\hat h)=\eta$.

We change the variables to $\sigma=\log h$ and $\tau=\log\hat h$. Problem
\eqref{eq:design} becomes
\begin{equation}\label{eq:design-log}
 \text{minimize }\ K_1e^{-\sigma}+K_2e^{-\tau}
 \quad\text{subject to}\quad G_0\le0,\quad G_1\le0,\quad G_2\le0,
\end{equation}
with $G_0=\Phi(e^\sigma)+\Psi(e^\tau)-\eta$, $G_1=\sigma-\tau$ and $G_2=\tau-\log\bar h$. By
\eqref{eq:PhiPsi}, the Hessian matrix of $G_0$ is positive definite on
the whole plane, and so is that of the cost function. The functions
$G_1$ and $G_2$ are linear. Hence the feasible set is convex, and the cost
function is strictly convex. Thus the minimizing pair $(h^*,\hat h^*)$ is unique.

In the original variables the active constraint gradients at this pair are
linearly independent. Indeed, \eqref{eq:hbar} rules out simultaneous
activity of the two mesh constraints, and the accuracy gradient
$(\Phi'(h^*),\Psi'(\hat h^*))$ has two positive components. It is nonzero
and linearly independent of either mesh gradient $(1,-1)$ or $(0,1)$.
By \cite[Theorem~4.2.13]{BSS}, the minimizing pair therefore admits
nonnegative multipliers satisfying \eqref{eq:KKT}. To see uniqueness,
let $\xi$ and $\zeta$ be two such multiplier vectors, let
$x^*=(h^*,\hat h^*)$, and set $I=\{i:g_i(x^*)=0\}$. Complementarity
gives $\xi_i=\zeta_i=0$ for $i\notin I$. Subtracting the stationarity
equations gives
\[
 \sum_{i\in I}(\xi_i-\zeta_i)\nabla g_i(x^*)=0.
\]
Linear independence implies $\xi_i=\zeta_i$ for every $i\in I$.

Conversely, if $(h,\hat h;\xi_0,\xi_1,\xi_2)$ is a KKT solution of
\eqref{eq:design}, then its image is a KKT solution of \eqref{eq:design-log}
under the map
\[
 (h,\hat h;\xi_0,\xi_1,\xi_2)
 \longmapsto
 (\log h,\log\hat h;\xi_0,h\xi_1,\hat h\xi_2).
\]
Indeed, the map preserves feasibility, nonnegativity of the multipliers
and complementarity. Since
$\xi_1(h-\hat h)=0$, we have $h\xi_1=\hat h\xi_1$. Multiplying the two
stationarity equations in \eqref{eq:KKT} by $h$ and $\hat h$, respectively,
therefore gives the stationarity equations for \eqref{eq:design-log}.
By \cite[Theorem~4.2.16(b)]{BSS}, the transformed pair minimizes
\eqref{eq:design-log}. Thus $(h,\hat h)=(h^*,\hat h^*)$, and the
preceding uniqueness argument determines the multipliers as well.
\end{proof}

\begin{proof}[Proof of Corollary~\ref{cor:design-structure}]
The functions $\Phi$ and $\Psi$ have positive derivatives, so they
increase strictly. By Theorem~\ref{prop:design},
$\Phi(h^*)+\Psi(\hat h^*)=\eta=\Phi(h_1)+\Psi(h_1)$. If $h^*>h_1$, then
$\hat h^*\ge h^*>h_1$, and $\Phi(h^*)+\Psi(\hat h^*)>\eta$. If
$\hat h^*<h_1$, then $h^*\le\hat h^*<h_1$, and
$\Phi(h^*)+\Psi(\hat h^*)<\eta$. Both contradict the equality. This
proves (a).

By \eqref{eq:hbar}, $\Phi(h_1)+\Psi(h_1)=\eta<\Phi(\bar h)+\Psi(\bar h)$,
so $h_1<\bar h$. At $h=\hat h=h_1$, the last condition in \eqref{eq:KKT}
gives $\xi_2=0$. Adding the first two equations gives
$\xi_0=(K_1+K_2)/(h_1^2(\Phi'(h_1)+\Psi'(h_1)))$, and then
$\xi_1=K_1/h_1^2-\xi_0\Phi'(h_1)=(K_1\Psi'(h_1)-K_2\Phi'(h_1))/(h_1^2(\Phi'(h_1)+\Psi'(h_1)))$.
So the pair $(h_1,h_1)$ satisfies \eqref{eq:KKT} if and only if
$K_2\Phi'(h_1)\le K_1\Psi'(h_1)$. By Theorem~\ref{prop:design}, this
holds if and only if $(h_1,h_1)$ is the minimizer. This proves (b).

If $h^*<\hat h^*<\bar h$, then $\xi_1=\xi_2=0$ by the complementarity
conditions, and the first two
equations of \eqref{eq:KKT} give
$K_1/((h^*)^2\Phi'(h^*))=\xi_0=K_2/((\hat h^*)^2\Psi'(\hat h^*))$. This is
the equality in (c).
With $x=(\hat h^*)^{1/6}$, the equality $\Phi(h^*)+\Psi(\hat h^*)=\eta$ gives
$h^*=t_0^{3/2}(\eta-\Psi(x^6))/C_1$, so $\Psi(x^6)<\eta$ and the formula
for the ratio holds. Inserting $h^*$ into the equality in (c) gives
$P(x)=0$, since
$(\hat h^*)^2\Psi'(\hat h^*)=\frac12C_2t_0^{-5/2}\varepsilon x^9+\frac16C_3\varepsilon^2x^7$.
Where $y>0$ and $\Psi(y^6)<\eta$, the equation $P(y)=0$ is equivalent to
$(K_1C_1/(K_2t_0^{3/2}))^{1/2}((y^6)^2\Psi'(y^6))^{1/2}+\Psi(y^6)=\eta$,
whose left side increases strictly in $y$. So $x$ is the only such zero.
This proves (c).

For (d), the equality in (c) reads
$K_1t_0^{3/2}/(C_1h^2)=K_2/(\hat h^2\Psi'(\hat h))$. It gives the formula
for the ratio and $h=G(\hat h)$ with
$G(y)=(K_1t_0^{3/2}y^2\Psi'(y)/(C_1K_2))^{1/2}$. Since $y^2\Psi'(y)$
increases strictly, so does $G$. For each $\eta$ in the interval of (d),
$\hat h^*$ solves
$C_1t_0^{-3/2}G(y)+\Psi(y)=\eta$, whose left side increases strictly in
$y$. Hence $\hat h^*$ increases strictly with respect to $\eta$, and so does
$h^*=G(\hat h^*)$. The exponents $1/2$ and $1/6$ in $\Psi$ are less than
$1$, so $\Psi'$ decreases strictly. Hence the ratio increases strictly
with respect to $\hat h^*$.
\end{proof}

\begin{remark}[marginal costs]\label{rem:marginal}
Let $s$ be the share of the error budget given to the linear part, so
that $h=\Phi^{-1}(s)$. Then
$K_1/(h^2\Phi'(h))=-\,\mathrm d(K_1/h)/\mathrm ds$ is the cost saved by the
last unit of error given to the fine mesh. If $h^*<\hat h^*<\bar h$, then
by Corollary~\ref{cor:design-structure}(c), at the optimum the last unit
of error saves the same cost on both meshes. The common value is the
multiplier $\xi_0$ of the accuracy constraint.
\end{remark}

The condition in Corollary~\ref{cor:design-structure}(b) means that the
expression for the ratio in (d), taken at $\hat h=h_1$, is at most $1$.
When it holds, the two-mesh computation costs
$2J_{\rm L}(h_1)$ more than one nonlinear run on the mesh $h_1$. Problem
\eqref{eq:design} does not contain the single run. The two-mesh
computation is the cheaper one if and only if
$J_0+K_1/h^*+K_2/\hat h^*<J_{\rm N}(h_1)$.

The next corollary concerns small data. The nonlinear effect $u-v$ is at
most of order $\varepsilon^2$ (Theorem~\ref{thm:linear}). We ask for an
error of at most $\rho\varepsilon^2$ with a fixed $\rho>0$. In
\eqref{eq:accuracy} this means $\eta=\rho\varepsilon$.

\begin{corollary}[small data]\label{cor:design-small}
Let $d_{\rm N}>d_{\rm L}$. Let $\Phi$ and $\Psi$ be defined as in
\eqref{eq:PhiPsi}. Fix $\rho>0$ and $\bar h$, and put
$\eta=\rho\varepsilon$. Let $J^*$ be the least serial cost of $\tilde u$
under \eqref{eq:design}, and let $h_{1,\varepsilon}$ solve
$\Phi(h)+\Psi(h)=\rho\varepsilon$. As $\varepsilon\to0$,
\[
 \frac{J^*}{J_{\rm N}(h_{1,\varepsilon})}\to\frac{d_{\rm L}}{d_{\rm N}}<1
 \qquad\text{and}\qquad \frac{\hat h^*}{h^*}\to\infty .
\]
\end{corollary}

\begin{proof}
For small $\varepsilon$ we have $\Phi(\bar h)>\rho\varepsilon$, so
Theorem~\ref{prop:design} applies, and $h_{1,\varepsilon}\le\bar h$.
Since $\Phi+\Psi$ is increasing, $h_{1,\varepsilon}$ is the coarsest
single mesh admitted by the bound. As $\varepsilon\to0$, we have
$h_{1,\varepsilon}\to0$ and $\Psi(h_{1,\varepsilon})=o(\varepsilon)$, so
$h_{1,\varepsilon}=\rho t_0^{3/2}C_1^{-1}\varepsilon\,(1+o(1))$. By
Corollary~\ref{cor:design-structure}(a), $h^*\le h_{1,\varepsilon}$, so
$J^*\ge J_0+K_1/h_{1,\varepsilon}+K_2/\hat h^*$. For an upper bound, take
$\hat h_\varepsilon=\varepsilon^{1/2}$. Then
$\Psi(\hat h_\varepsilon)=O(\varepsilon^{5/4})=o(\varepsilon)$. Let
$h_\varepsilon$ solve $\Phi(h_\varepsilon)=\rho\varepsilon-\Psi(\hat h_\varepsilon)$.
Then $h_\varepsilon=h_{1,\varepsilon}(1+o(1))$ and
$h_\varepsilon\le\hat h_\varepsilon\le\bar h$ for small $\varepsilon$. The
pair is a feasible solution, and
\[
 J^*\le J_0+\frac{K_1}{h_\varepsilon}+\frac{K_2}{\hat h_\varepsilon}
 =J_0+\frac{K_1}{h_{1,\varepsilon}}\,(1+o(1))+O(\varepsilon^{-1/2}).
\]
Since $J_{\rm N}(h_{1,\varepsilon})=J_0'+K/h_{1,\varepsilon}$ and
$K/h_{1,\varepsilon}$ is of order $\varepsilon^{-1}$, the two bounds give
$J^*/J_{\rm N}(h_{1,\varepsilon})\to K_1/K=d_{\rm L}/d_{\rm N}$. They also
give $K_2/\hat h^*=o(\varepsilon^{-1})$. Hence
$\hat h^*/\varepsilon\to\infty$, while $h^*\le h_{1,\varepsilon}=O(\varepsilon)$.
\end{proof}

\begin{remark}[parallel computation]\label{rem:parallel}
The three parts of $\tilde u$ do not depend on each other. Run on three
processors, they take the wall-clock time
$J_\parallel=\max\{J_{\rm L}(h),J_{\rm L}(\hat h),J_{\rm N}(\hat h)\}$.
Every feasible solution has $h\le h_1$, since
$\Phi(h)+\Psi(h)\le\Phi(h)+\Psi(\hat h)\le\eta$. So
$J_\parallel\ge J_{\rm L}(h_1)$. In the setting of
Corollary~\ref{cor:design-small}, the pair
$(h_\varepsilon,\hat h_\varepsilon)$ of its proof gives
$J_\parallel=J_{\rm L}(h_\varepsilon)$ for small $\varepsilon$, and
$h_\varepsilon=h_{1,\varepsilon}(1+o(1))$. Hence
$\min J_\parallel/J_{\rm L}(h_{1,\varepsilon})\to1$ as $\varepsilon\to0$:
the least wall-clock time is asymptotically that of the linear simulation
on the fine mesh alone.
\end{remark}

\subsection{\texorpdfstring{Data in $D(A)$}{Data in D(A)}}\label{sec:design-smooth}
For smooth data both error terms are of second order, and the bounds hold
down to $t=0$. Let $\nu$ be as in Section~\ref{sec:design-constraint},
and let $\E_\nu$ be the energy norm \eqref{eq:Enorm}. On $D(A)$ the graph
norm of $A$ is equivalent to the $H^3$ norm. By
Theorem~\ref{thm:DA}(i),
\begin{equation}\label{eq:smoothlinear}
 \E_\nu(v-v_h)\le C_4h^2\norm{Au_0}_H\qquad(u_0\in D(A)).
\end{equation}
By Theorem~\ref{thm:DA}(ii), the same bound holds for $u-u_h$ when
$\norm{u_0}_H$ is small. Hence $\E_\nu(\Cd_{\hat h})\le C\hat h^2\norm{Au_0}_H$.
This bound has no factor of the size of the datum.
Theorem~\ref{thm:DA}(iii) supplies it: there is $C_5$, depending only on
$L$ and $\nu$, such that
\begin{equation}\label{eq:smoothordering}
 \E_\nu(\Cd_{\hat h})\le C_5\hat h^2\norm{Au_0}_H^{3/2}\norm{u_0}_H^{1/2}
 \qquad(u_0\in D(A)\text{ real},\ \norm{u_0}_H\le r).
\end{equation}

Let $M>0$, and consider data with $\norm{Au_0}_H\le M\norm{u_0}_H$. Then
the right sides of \eqref{eq:smoothlinear} and \eqref{eq:smoothordering}
are at most $C_4Mh^2\norm{u_0}_H$ and $C_5M^{3/2}\hat h^2\norm{u_0}_H^2$.
On the same mesh, the bound for the ordering error has the extra factor
$\norm{u_0}_H$. For small data, the coarse mesh can therefore be much
coarser than the fine mesh (Theorem~\ref{cor:design-smooth}(c)). The constants $C_4,C_5$ depend only on $L$ and $\nu$. They do not
depend on $h$, $\hat h$, $t$, $M$, $\varepsilon$, $\eta$ and $u_0$.
Now let $\norm{u_0}_H\le\varepsilon\le r$ and $h\le\hat h$. By
\eqref{eq:smoothlinear}, \eqref{eq:smoothordering} and
$u-\tilde u=(v-v_h)+\Cd_{\hat h}$, the bound
\begin{equation}\label{eq:smoothPhiPsi}
 \Phi(h)+\Psi(\hat h)\le\eta,\qquad \Phi(h)=ah^2,\quad\Psi(\hat h)=c\hat h^2,
 \quad a=C_4M,\quad c=C_5M^{3/2}\varepsilon,
\end{equation}
implies $\norm{(u-\tilde u)(t)}_H\le\eta\norm{u_0}_H$ for all $t\ge0$. No
$t_0$ is needed. Now $\Phi$ and $\Psi$ are convex, and so is problem
\eqref{eq:design}. Its solution is explicit. Put
\[
 k=\frac{K_2}{K_1}=1+\frac{d_{\rm N}}{d_{\rm L}},\qquad
 \gamma=\Big(k\,\frac ac\Big)^{1/3},\qquad
 R(\gamma)=\frac{d_{\rm L}}{d_{\rm N}}\cdot\frac{(\gamma+k)^{3/2}}{(\gamma^3+k)^{1/2}} .
\]

\begin{theorem}[Design for data in $D(A)$]\label{cor:design-smooth}
Let $d_{\rm N}>d_{\rm L}$. Let $\Phi(h)=ah^2$ and $\Psi(\hat h)=c\hat h^2$
with $a,c>0$, and let $h_1=(\eta/(a+c))^{1/2}<\bar h$.
\begin{enumerate}
\item[(a)] If $\gamma\le1$, then $h^*=\hat h^*=h_1$ for every $\eta$.
\item[(b)] Let $\gamma>1$, and put $h_\gamma=(\eta/(a+c\gamma^2))^{1/2}$. If
$\gamma h_\gamma\le\bar h$, then $h^*=h_\gamma$ and $\hat h^*=\gamma h_\gamma$. Thus the
mesh ratio $\hat h^*/h^*=\gamma$ does not depend on $\eta$. Moreover
\[
 \frac{K_1/h^*+K_2/\hat h^*}{K/h_1}=R(\gamma),
\]
and $R$ decreases strictly on $(1,\infty)$, from $(K_1+K_2)/K$ to
$d_{\rm L}/d_{\rm N}$.
\item[(c)] Let $a$ and $c$ be as in \eqref{eq:smoothPhiPsi}, and let
$\eta=\rho\varepsilon$ with a fixed $\rho>0$. As $\varepsilon\to0$, the
case (b) applies for small $\varepsilon$,
$\hat h^*/h^*=\gamma=(kC_4/C_5)^{1/3}M^{-1/6}\varepsilon^{-1/3}$, and
$J^*/J_{\rm N}(h_1)\to d_{\rm L}/d_{\rm N}<1$.
\end{enumerate}
\end{theorem}

\begin{proof}
With these $\Phi$ and $\Psi$, problem \eqref{eq:design} is convex: the
objective is strictly convex on $(0,\infty)^2$, and the constraints are
convex. So it has at most one minimizer, and a feasible solution that satisfies \eqref{eq:KKT} is the
minimizer~\cite[Theorem~4.2.16(b)]{BSS}. For (a), take
$h=\hat h=h_1$, $\xi_0=(K_1+K_2)/(2(a+c)h_1^3)$,
$\xi_1=K_1/h_1^2-2a\xi_0h_1$ and $\xi_2=0$. All equations in
\eqref{eq:KKT} hold, and $\xi_1=(K_1c-K_2a)/((a+c)h_1^2)$ is nonnegative
if and only if $ka\le c$, that is, $\gamma\le1$. This proves (a). For (b),
the pair $(h_\gamma,\gamma h_\gamma)$ is a feasible solution and
$\Phi(h_\gamma)+\Psi(\gamma h_\gamma)=\eta$. Since $c\gamma^3=ka$,
\eqref{eq:KKT} holds with $\xi_0=K_1/(2ah_\gamma^3)$ and
$\xi_1=\xi_2=0$. So the pair is the minimizer. From $c=ka/\gamma^3$ we
get $a+c\gamma^2=a(\gamma+k)/\gamma$ and $a+c=a(\gamma^3+k)/\gamma^3$.
Hence
\[
 \frac{K_1/h_\gamma+K_2/(\gamma h_\gamma)}{K/h_1}
 =\frac{K_1}K\Big(1+\frac k\gamma\Big)\Big(\frac{a+c\gamma^2}{a+c}\Big)^{1/2}
 =R(\gamma),
\]
since $K_1/K=d_{\rm L}/d_{\rm N}$. Moreover
$\frac{d}{d\gamma}\log R(\gamma)=\frac32k(1-\gamma^2)/((\gamma+k)(\gamma^3+k))<0$
for $\gamma>1$, $R(1)=(K_1+K_2)/K$ and $R(\gamma)\to K_1/K$ as
$\gamma\to\infty$. For (c), $c\to0$ while $a$ is fixed. Then
$\gamma=(kC_4/C_5)^{1/3}M^{-1/6}\varepsilon^{-1/3}$ and
$h_\gamma\le(\eta/a)^{1/2}=(\rho\varepsilon/a)^{1/2}$. Hence
$\gamma h_\gamma\le C\varepsilon^{1/6}\to0$, and also $h_1\to0$. Hence (b)
applies for small $\varepsilon$, and $R(\gamma)\to d_{\rm L}/d_{\rm N}$.
The mesh-dependent costs tend to infinity, so the fixed costs $J_0$ and
$J_0'$ do not change the limit.
\end{proof}

In case (b), and apart from the costs that do not depend on the meshes,
the two-mesh computation is the cheaper one if and only if
$\gamma>\gamma_c$, where $\gamma_c>1$ solves $R(\gamma_c)=1$, that is,
$(\gamma_c+k)^3=(d_{\rm N}/d_{\rm L})^2(\gamma_c^3+k)$. For
$d_{\rm N}/d_{\rm L}=3.8$, $\gamma_c\approx3.07$. So, while the coarse
mesh stays below $\bar h$, the mesh ratio for data in $D(A)$ does not
depend on the accuracy. For fixed error constants $C_4,C_5$, it is
determined by $M$, $\varepsilon$ and $d_{\rm N}/d_{\rm L}$. When
$\gamma h_\gamma>\bar h$, the coarse mesh is $\bar h$ and the fine mesh is
$((\eta-c\bar h^2)/a)^{1/2}$, so the ratio decreases as $\eta$ grows. This occurs for
fixed $0<\eta\le a\bar h^2$ and sufficiently small $\varepsilon$.
The limiting ratio of the mesh-dependent costs then exceeds
$d_{\rm L}/d_{\rm N}$. For fixed $\eta>a\bar h^2$, the coefficient
$c=C_5M^{3/2}\varepsilon$ tends to zero, so
$(a+c)\bar h^2\le\eta$ for all sufficiently small $\varepsilon$.
Both optimal mesh sizes then equal $\bar h$.
For data in a ball of $H$, the optimal mesh ratio increases with $\eta$
while both mesh constraints are inactive
(Corollary~\ref{cor:design-structure}(d)), and one nonlinear run is cheaper when
$\eta$ is at most a small constant times $\varepsilon^2$
(Remark~\ref{rem:regimes}).

\subsection{Regimes and constants}\label{sec:design-regimes}

\begin{remark}[three regimes]\label{rem:regimes}
High accuracy and a coarse correction mesh are not compatible. The
admitted error for data of size $\varepsilon$ is $\eta\varepsilon$, and the
bounds select one of three computations. Let $C_0$ be the constant
in $\norm{(u-v)(t)}_H\le C_0\norm{u_0}_H^2$ (Theorem~\ref{thm:linear}).
\begin{enumerate}
\item[(i)] If $\eta>C_0\varepsilon$, the linear simulation alone suffices.
Since $u-v_h=(u-v)+(v-v_h)$, the bound gives \eqref{eq:accuracy} for
$v_h$ in place of $\tilde u$ as soon as $\Phi(h)\le\eta-C_0\varepsilon$. No
correction is needed.
\item[(ii)] If $K_2\Phi'(h_1)\le K_1\Psi'(h_1)$, the optimum of
\eqref{eq:design} lies on $h=\hat h$, and one nonlinear run on the mesh
$h_1$ is cheaper than the two-mesh computation. For \eqref{eq:PhiPsi} we
have $\Psi'(h)\ge\frac12C_2t_0^{-5/2}\varepsilon h^{-1/2}$ and
$h_1\le\eta t_0^{3/2}/C_1$. With $a=C_1t_0^{-3/2}$, $b=C_2t_0^{-5/2}$
and $k=K_2/K_1$, the condition therefore holds whenever
$\eta\le b^2\varepsilon^2/(4ak^2)$, that is, whenever the admitted error is
at most a small constant times $\varepsilon^3$. For
\eqref{eq:smoothPhiPsi} the condition is $\gamma\le1$, whatever $\eta$
(Theorem~\ref{cor:design-smooth}(a)).
\item[(iii)] In between, the two-mesh computation with the optimal pair
of Theorem~\ref{prop:design} is the candidate. It is cheaper than one
nonlinear run when $J_0+K_1/h^*+K_2/\hat h^*<J_{\rm N}(h_1)$. For
$\eta=\rho\varepsilon$ with $\rho<C_0$, case (i) does not apply, and
Corollary~\ref{cor:design-small} and Theorem~\ref{cor:design-smooth}(c) show
that this holds for small data.
\end{enumerate}
Thus, for data in a ball of $H$, the optimal correction mesh may be
coarser than the fine mesh only when the admitted error exceeds a
constant times $\varepsilon^3$. For data in $D(A)$, while the coarse mesh
stays below $\bar h$, the accuracy does not enter the optimal mesh ratio.
\end{remark}

\begin{remark}[constants]\label{rem:design-constants}
The constants $C_1,\dots,C_5$ are not explicit. The design above gives the
structure of the optimum, not its numerical values. For data in $D(A)$,
in case (b) of Theorem~\ref{cor:design-smooth}, the optimal mesh ratio
depends on the constants only through $a/c$. The
optimal meshes depend on the cost parameters only through
$d_{\rm N}/d_{\rm L}$. The comparison with one nonlinear run also involves
the fixed works per step, through $J_0-J_0'=N_t(2c_{\rm L}-2d_{\rm L})$.
Section~\ref{sec:comp} fixes the constants on coarse meshes and checks
the resulting designs.
\end{remark}

\section{Experimental design and results}\label{sec:comp}
\setcounter{equation}{0}

In this section, we test the designs of Section~\ref{sec:design} for the
W-scheme with $L=1$. The computations illustrate the results. They do not
prove them.

\begin{figure}[H]
\centering
\includegraphics[width=\textwidth]{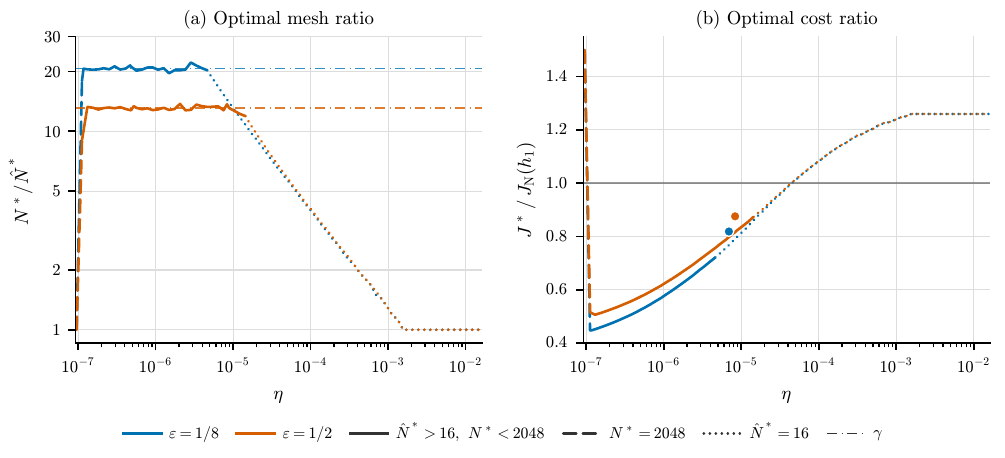}
\caption{Designs for D1. (a) Integer optimal mesh ratio $N^*/\hat N^*$ and
the ratio $\gamma$ of Theorem~\ref{cor:design-smooth} (thin
dash-dotted lines). (b) Least model cost $J^*$ of the two-mesh computation
divided by the cost $J_{\rm N}(h_1)$ of one nonlinear run on $N_1$. Lines:
computed from the bounds \eqref{eq:smoothPhiPsi} with the constants of
Table~\ref{tab:smooth-const} and the cost model; dots: measured at the
designed pairs of Table~\ref{tab:smooth-check}.}\label{fig:smooth}
\end{figure}

\paragraph{Three stages.} Each experiment has three stages.
\begin{enumerate}
\item \emph{Calibration.} The constants of the error bounds are fixed
from runs on the coarse meshes $N,\hat N\in\{16,32\}$. The parameters of
the cost model are fixed from short runs with $N\le512$. These runs are
not used again.
\item \emph{Design.} For each datum, amplitude $\varepsilon$ and accuracy
$\eta$, we compute from the bounds and the cost model the optimal pair
$(N^*,\hat N^*)$, and the coarsest mesh $N_1$ on which the same bounds
admit one nonlinear run. This stage solves an optimization problem. It
requires no time integration of the KdV equation. Its results are the lines in Figure~\ref{fig:smooth}.
\item \emph{Check.} At designed pairs chosen by a rule fixed in advance,
we run the scheme and measure errors and times.
\end{enumerate}
Before any check, the constants, the designed pairs, the predicted costs
and the predictions below were written to files. Their hashes and times
are recorded in the ancillary files. Failed predictions are reported.

\paragraph{The runs at a designed pair.} At a pair $(N,\hat N)$ with single
mesh $N_1$, we use four runs of the W-scheme, each from $t=0$ to $t=1$:
the linear scheme on $N$ and on $\hat N$, and the nonlinear scheme on
$\hat N$ and on $N_1$. The first three runs give $\tilde u$. The last
gives $u_{h_1}$. All runs use the five-stage Radau~IIA method, with time
grids generated by one rule. Identical runs are shared between pairs,
and each run is timed separately.

\paragraph{Measured quantities.}
\begin{itemize}
\item \emph{Error.} For $w=\tilde u$ or $w=u_{h_1}$, the error is the
maximum of $e^{\nu t}\norm{(u-w)(t)}_H$, $\nu=3\pi^2/8$, over the output
times, with $u$ replaced by a reference solution on $N=4096$. The output times
cover $[0,1]$ for data in $D(A)$, and $[t_0,1]$ for data in a ball of
$H$. We divide the error by the admitted error $\eta\varepsilon$. A
quotient of at most $1$ meets the requirement at the output times. At $t=0$
the datum itself is used. The weight $e^{\nu t}\ge1$ makes the test at
least as strict as the unweighted test against the same reference.
\item \emph{Cost ratio.} The serial time of $\tilde u$ is the sum of the
times of its three runs. The cost ratio is this time divided by the time
of the run on $N_1$. The predicted ratio is the same quotient, with times
from the cost model and the numbers of time steps fixed in advance by
the rule for the time grids. The
time grid depends on the datum, and on the mesh only through an initial
layer. A ratio below $1$ means that the two-mesh
computation is cheaper.
\end{itemize}

\paragraph{Rules fixed in advance.}
\begin{itemize}
\item \emph{Resolution.} For data in $D(A)$ the empirical reference
uncertainty for the full solution is $U=S+T_{2048}+T_{4096}+F$.
Here $S$ is the difference under spatial refinement from $2048$ to $4096$
elements at time-grid level~3. The differences $T_{2048}$ and $T_{4096}$
are obtained by refining the time grid from level~2 to level~3 on the
two meshes separately. Each difference is measured by the maximum of
the exponentially weighted $H$-norm over positive output times. The
roundoff floor is
\[
 F=64\,\varepsilon_{\rm mach}
 \max\{\varepsilon,\max_t e^{\nu t}\norm{u_{\rm ref}(t)}_H\},
\]
where the maximum includes $t=0$. An error is resolved if it is at
least $10U$.
\item \emph{Selection.} For each smooth datum and amplitude we tried
pairs with $N$ near $256$, $1024$ and $2048$. We kept a pair if
$\eta\varepsilon\ge100U$ and either $S\le0.8S_{\rm prev}$ or $S\le F$.
Here $S_{\rm prev}$ is the corresponding spatial difference from
$1024$ to $2048$ elements.
\item \emph{Predictions.} For data in $D(A)$: (P1) both errors are at most
$\eta\varepsilon$; (P2) when both the integer and the continuous
optimal pairs are away from the mesh limits, the integer optimal ratio
$N^*/\hat N^*$ lies within $15\%$ of $\gamma$; (P3) the measured
cost ratio lies within $25\%$ of the predicted one; (P4) $\gamma$
decreases as $\varepsilon$ grows.
\end{itemize}

\paragraph{Cost parameters.}
We fitted the time per step, $c+dn$ with $n=2N-1$ unknowns, to the
measured times of linear and nonlinear runs on $N=16,\dots,512$, timed
on one thread of an eight-core Intel machine. The times include the
factorizations of the banded matrices and exclude their assembly.
Table~\ref{tab:cost} gives the result.

\begin{table}[H]
\centering\small
\caption{Cost parameters: seconds per time step, $c+dn$.}\label{tab:cost}
\begin{tabular}{@{}lcc@{}}
\toprule
run & $c$ & $d$\\
\midrule
linear & $1.54\times10^{-4}$ & $8.02\times10^{-7}$\\
nonlinear & $1.28\times10^{-3}$ & $3.05\times10^{-6}$\\
\bottomrule
\end{tabular}
\end{table}

Thus $d_{\rm N}/d_{\rm L}=3.80$, $k=4.80$ and $\gamma_c\approx3.07$. A
nonlinear step needs about two fixed-point iterations. Its fixed work
$c_{\rm N}$ equals the work for about $420$ unknowns. So the fraction
$d_{\rm L}/d_{\rm N}=0.26$ of Corollary~\ref{cor:design-small}
and Theorem~\ref{cor:design-smooth} is approached only on meshes with many more
unknowns. The predicted times of runs on $N=1024$ and $2048$, written
before these runs, deviated from the measured ones by $4\%$ to $37\%$ of the measured
times. The
parameters depend on the machine and on the implementation.

\paragraph{Data in $D(A)$.}
We use four smooth data with $\norm f_H=1$: D1, a bump supported in
$[0.25,0.75]$; D5, a narrower bump supported in $[0.4,0.6]$; and, held
out, D1$'$, a bump supported in $[0.2,0.65]$, and D6, the bump of D1
times $\sin6\pi x$. With $\mathcal B(s)=\exp(-1/s-1/(1-s))$ for $0<s<1$
and $\mathcal B(s)=0$ otherwise, they are $\mathcal B((x-0.25)/0.5)$,
$\mathcal B((x-0.4)/0.2)$, $\mathcal B((x-0.2)/0.45)$ and
$\mathcal B((x-0.25)/0.5)\sin6\pi x$, each divided by its norm in $H$.
Their values of $M=\norm{Af}_H$ are
$3.5\times10^3$, $5.5\times10^4$, $4.8\times10^3$ and $1.9\times10^4$. The
amplitudes are $\varepsilon=1/32$, $1/8$ and $1/2$. The constants $a$ and
$c$ of \eqref{eq:smoothPhiPsi} were fixed on the coarse meshes
$N,\hat N\in\{16,32\}$, for D1 and D5, as the largest values of
$N^2\times$(linear error) and of $\hat N^2\times$(ordering
error)$/\varepsilon$. No safety factor was used. D1$'$ uses the constants
of D1, and D6 those of D5. Table~\ref{tab:smooth-const} lists them.

\begin{table}[H]
\centering\small
\caption{Constants of \eqref{eq:smoothPhiPsi} for unit data, fixed on
$N,\hat N\in\{16,32\}$, and the mesh ratio $\gamma$ of
Theorem~\ref{cor:design-smooth}.}\label{tab:smooth-const}
\begin{tabular}{@{}lcccccccc@{}}
\toprule
 & & \multicolumn{3}{c}{$c/\varepsilon$} & \multicolumn{3}{c}{$\gamma$}\\
\cmidrule(lr){3-5}\cmidrule(l){6-8}
datum & $a$ & $\varepsilon=1/32$ & $1/8$ & $1/2$ & $\varepsilon=1/32$ & $1/8$ & $1/2$\\
\midrule
D1 & $0.399$ & $1.719\times10^{-3}$ & $1.719\times10^{-3}$ & $1.719\times10^{-3}$ & $32.9$ & $20.7$ & $13.1$\\
D5 & $7.39$ & $2.322\times10^{-2}$ & $2.322\times10^{-2}$ & $2.321\times10^{-2}$ & $36.6$ & $23.0$ & $14.5$\\
\bottomrule
\end{tabular}
\end{table}

The quotient $c/\varepsilon$ agrees to three digits for the three
amplitudes. This is consistent with the factor $\varepsilon$ of
Theorem~\ref{thm:DA}(iii): on these data the calibrated ordering error
scales quadratically with the amplitude. The ratio $\gamma$ falls by the
factor $4^{1/3}\approx1.59$ from one amplitude to the next, and it exceeds
$\gamma_c$ in all cases. On the finer meshes $N=64,\dots,512$, the
measured errors of D1 and D5 stayed below the frozen models $ah^2$ and
$c\hat h^2$. For D1$'$, the measured linear and ordering errors were below their
respective model values on all four finer meshes ($N=64,128,256,512$).
The largest ratios over all tested meshes were $1.53$ and $1.39$, respectively,
both at $N=16$. For D6, both normalized errors were below $30\%$
of their respective model values.

\paragraph{Designs.}
For each datum and amplitude we minimized the cost model over the integer
pairs $16\le\hat N\le N\le2048$ under $ah^2+c\hat h^2\le\eta$, and we
compared with one nonlinear run on the coarsest admitted mesh $N_1$.
Figure~\ref{fig:smooth} shows D1. Over the twelve curves, $70{,}294$ integer
optima have both the integer and the continuous optimal pairs away from
the limits $\hat N=16$ and $N=2048$. Their ratios $N/\hat N$ lie within
$15\%$ of $\gamma$, with a largest deviation of $11.6\%$. Near the limits
the deviation can be larger. Once $\hat N=16$, the ratio decreases
as $\eta$ grows, as noted after Theorem~\ref{cor:design-smooth}. For
$\varepsilon=1/8$ the model cost ratio rises from $0.45$ near $N=2048$ to
$1$ at $\eta\approx4.5\times10^{-5}$, where $N\approx96$. For larger
$\eta$ one nonlinear run is cheaper, because of the fixed work per step.

\paragraph{Checks.}
The selection rule kept $18$ pairs: the twelve pairs with $N$ near $256$
of Table~\ref{tab:smooth-check}, and six pairs with $N$ near $1024$ for D5
and D6. No pair with $N$ near $2048$ qualified. (P1) holds at all
resolved pairs. (P2) and (P4) hold on all twelve curves. Among the
resolved pairs, (P3) fails at the three pairs whose measured cost ratio
exceeds $1$. At the unresolved pairs, (P1) and (P3) are not judged.

\begin{table}[H]
\centering\small
\caption{Designed pairs: errors of the two-mesh sum $\tilde u$ and of one
nonlinear run on $N_1$, divided by the admitted error $\eta\varepsilon$,
and the ratio of their serial times, predicted and measured; pairs with
$N$ near $256$.}\label{tab:smooth-check}
\begin{tabular}{@{}llccccc@{}}
\toprule
datum & $\varepsilon$ & $N/\hat N/N_1$ & error $\tilde u$ & error $u_{h_1}$ & predicted & measured\\
\midrule
D1 & $1/32$ & $256/16/252$ & $0.137$ & $0.143$ & $0.749$ & $0.758$\\
D1 & $1/8$ & $256/16/241$ & $0.126$ & $0.142$ & $0.768$ & $0.819$\\
D1 & $1/2$ & $260/19/220$ & $0.234$ & $0.141$ & $0.816$ & $0.876$\\
D1$'$ & $1/32$ & $256/16/252$ & $0.168$ & $0.174$ & $0.749$ & $0.820$\\
D1$'$ & $1/8$ & $256/16/241$ & $0.174$ & $0.173$ & $0.768$ & $0.840$\\
D1$'$ & $1/2$ & $260/19/220$ & $0.313$ & $0.173$ & $0.816$ & $0.867$\\
D5 & $1/32$ & $256/16/253$ & $0.053$ & $0.055$ & $0.748$ & $0.835$\\
D5 & $1/8$ & $256/16/245$ & $0.096$ & $0.057$ & $0.761$ & $1.508$\\
D5 & $1/2$ & $255/18/223$ & $0.178$ & $0.063$ & $0.806$ & $1.017$\\
D6 & $1/32$ & $256/16/253$ & $0.022$ & $0.022$ & $0.748$ & $0.880$\\
D6 & $1/8$ & $256/16/245$ & $0.027$ & $0.023$ & $0.761$ & $1.060$\\
D6 & $1/2$ & $255/18/223$ & $0.059$ & $0.025$ & $0.805$ & $0.969$\\
\bottomrule
\end{tabular}
\end{table}

All errors in Table~\ref{tab:smooth-check} are resolved, and they are
below the admitted ones by factors from $3$ to $45$. At the six pairs near
$N=1024$ the errors were below the admitted ones as well, but at each of
these pairs at least one error was not resolved. On meshes near $N=256$
the two-mesh computation
took $0.76$ to $0.88$ of the time of one nonlinear run for D1 and D1$'$,
against $0.75$ to $0.82$ predicted. For D5 and D6 the measured ratio was
$0.84$ to $1.51$, and at three pairs two meshes cost more than one. At the
pair with ratio $1.51$ the single run took $28$ seconds, against $41$ to
$44$ seconds at the other two amplitudes. We did not find the cause. Near
$N=1024$ the measured ratios were $0.50$ to $0.77$. For D6 they exceeded
the predicted ones by $37\%$ to $60\%$. Over all $18$ pairs, the measured
ratio differed from the predicted one by more than $25\%$ at six pairs. The measured absolute
times differed from the model by up to a factor of about two; the ratios
were closer. For D5, $\varepsilon=1/8$, $(N,\hat N)=(1022,45)$, the three parts
of $\tilde u$ run on three cores took $31$ seconds of wall-clock time,
including the start of the processes and the return of the results. This
is $0.39$ of the integration time of one nonlinear run on $N_1=932$.

On these meshes the savings are moderate, because of the fixed work per
step. On finer meshes the model cost ratio approaches $R(\gamma)$ of
Theorem~\ref{cor:design-smooth}, for example $R(20.7)=0.36$ for D1,
$\varepsilon=1/8$. The reference solution limits the checks to
$N\le1024$.

\paragraph{Data in a ball of $H$.}
For the bounds \eqref{eq:PhiPsi} we used the constants fixed
in [arXiv:2611.12585v2] on $N=16,32$, for D1 and D2 (a
function with
a jump), with $t_0=0.01,0.1$ and $\varepsilon=1/32,1/8$. With these
constants, $\Psi$ on the coarsest mesh $\hat N=16$ is smaller than $\Phi$
on the finest mesh $N=2048$ by two to four orders of magnitude. On each
computed design curve, the optimal fine mesh is the finest one, $N=2048$,
and the optimal coarse mesh reaches $\hat N=16$ at a tolerance less than
$0.5\%$ above the least feasible $\eta$. In these computations
$\Phi$ kept the minimum of \eqref{eq:sn-linear}, $\Phi(h)=C_1\min\{1,ht_0^{-3/2}\}$.
Since $\eta<C_1$ in all $30$ tested designs, this gives the same feasible
pairs as \eqref{eq:PhiPsi}. The $30$ tested designs are twelve for each of D1 and D2 and three for
each of two held-out data. All of them met the accuracy at the output
times, relative to the reference solution. The measured errors were at most $3.5\times10^{-4}$ of
the admitted error. The bound overestimates mainly the linear error. Its
factor $\min\{1,ht_0^{-3/2}\}$ equals $1$ for $N\le1000$ at $t_0=0.01$ and
decays only like $h$, whereas the computed linear errors of these data
decay about like $h^2$. At $(N,\hat N)=(2048,16)$ the measured cost ratio
was $0.44$ to $0.63$. For data in $D(A)$ the linear bound is of order
$h^2$, and the design is not degenerate.

\section*{Declarations}

\paragraph{Use of generative AI tools.}
During the preparation of this work, the author used Claude Opus~5.5 and
Claude Fable~5.2 (Anthropic), and GPT-6~Astra and GPT-6.2~Astra (OpenAI).
The author and the four models worked as a research team. The author set
every key goal and direction, discussed the plans of implementation with
the main thread, run by Claude Opus~5.5, and assigned tasks to other
threads, which carried them out in parallel. The tasks included drafts of
the text and of the proofs, literature search, the programs of the
computations, and a Lean~4 check of Theorems~\ref{thm:linear},
\ref{thm:ordering} and~\ref{thm:DA} under certain premises. GPT-6.2~Astra
and Claude Fable~5.2 cross-checked the drafts. On this basis, the author
checked and revised the draft line by line. At this stage Claude
Opus~5.5 served as a writing secretary. AI tools can produce incorrect
arguments and citations, and a human reader may also overlook some of
them. The workflow therefore contains two safeguards: the cross-checks by
two models and the Lean~4 check. The author takes full responsibility for
the content of the
article, except that the programs behind the numerical results were not
checked line by line by the author. The author nevertheless has confidence in the numerical output
because the computations are reproducible and the two models cross-checked the code.

\paragraph{Data availability.}
All programs, outputs and frozen records of the computations, and the
Lean~4 files, are in the ancillary files of the arXiv version of this
paper.

\end{document}